\documentstyle[12pt,amscd,amssymb]{amsart}

\def\DEBUG{n}

\newtheorem{th}{Theorem}[section] 
\newtheorem{lem}[th]{Lemma} 
 \newtheorem{prp}[th]{Proposition}

\theoremstyle{definition} 
\newtheorem{dfn}[th]{Definition} 
 \newtheorem{conj}[th]{Conjecture}
 \newtheorem{question}[th]{Question}

\theoremstyle{remark} 
\newtheorem{rem}[th]{Remark}

\newcommand{\rar}{\rightarrow} \newcommand{\dar}{\downarrow}
 \newcommand{\das}{\dashrightarrow}
\newcommand{\bfp}{{\Bbb{P}}} 
\newcommand{\bfq}{{\Bbb{Q}}} 
\newcommand{\bfz}{{\Bbb{Z}}} \newcommand{\bfg}{{\Bbb{G}}}
\newcommand{\bfr}{{\Bbb{R}}} \newcommand{\bfo}{{\Bbb{O}}}
\newcommand{\bfa}{{\Bbb{A}}} \newcommand{\bfn}{{\Bbb{N}}}
\newcommand{\co}{{\cal{O}}}
 
 \newcommand{\setmin}{\,{^{_\setminus}}\,}

\def\cO{{\cal O}}

\if\DEBUG y

\else

\fi

\newcommand{\Hom}{{\operatorname{Hom }}}
\newcommand{\spec}{{\operatorname{Spec\ }}}
\newcommand{\sym}{{\operatorname{Sym}}}

\newcommand{\Aut}{{\operatorname{Aut}}}

\newcommand{\reld}{\operatorname{rel.dim }}
 \newcommand{\red}{{\mbox{\small
red}}}

\newcommand{\lcm}{{\operatorname{lcm}}}
\newcommand{\dimension}{{\operatorname{dim}}}
\newcommand{\Star}{{\operatorname{Star}}}
\newcommand{\rk}{{\operatorname{Rank}}}

\begin{document}

\noindent%

\title[Weak semistable reduction in characteristic 0]{Weak
semistable reduction in characteristic 0\\ Preliminary Version}
\author[Dan Abramovich]{Dan Abramovich$^\dag$}
\thanks{$^\dag$ Partially supported by NSF grant DMS-9700520
and by an Alfred P. Sloan research fellowship.}
\author{Kalle Karu}
\address{Department of Mathematics\\ Boston University\\
111 Cummington\\ Boston, MA 02215, USA}
\email{abrmovic@@math.bu.edu}
\email{kllkr@@math.bu.edu}

\maketitle

\vspace{-.15in}

\addtocounter{section}{-1}
\section{INTRODUCTION}
Regretfully, we work over an algebraically closed field $k$ of
characteristic 0.
\subsection{The problem} Roughly speaking, the semistable reduction problem we
address here asks for the following:
\begin{quote}
Let $X\rar B$ be a surjective morphism of complex projective
varieties with geometrically integral generic fiber. Find a generically finite
proper surjective morphism (that is, an alteration) $B_1\rar B$, and a proper
birational morphism (that is, a modification) $Y\rar X\times_B B_1$, such that
the morphism $Y\rar B_1$ is nice.
\end{quote}
Of course, one needs to decide what a ``nice morphism'' means.

The question was posed, among other places, in the introduction of
\cite{te}, p. vii. It can be viewed as a natural extension of
Hironaka's theorem on resolution of singularities, which is in a sense
``the general fiber'' of semistable reduction.

\subsection{Brief history} The case { $\dim B=1, \dim
X=2$ is very old, see \cite{aw}. When $\dim B=1,$ semistable reduction
was obtained in \cite{te}, in the best possible sense: $Y$ is
nonsingular, and all the fibers are reduced, strict divisors of normal
crossings.

Using a result of Kawamata on ramified covers (see \cite{kawamata},
theorem 17), one can obtain semistable reduction ``in codimension 1''
over a base of arbitrary dimension.  Below, we will refer to the
result of Kawamata as {\bf ``Kawamata's trick''.} We will discuss it
in detail in section \ref{reduced-fibers}.

The case where $\dim X = \dim B + 1$ has recently been proven by de
Jong \cite{dj}. Here one shows that any family of curves can be made
into a family of nodal curves, which are indeed as ``nice'' as one may
expect.

Using recent difficult results of Alexeev, Koll\'ar and
Shepherd-Barron (see \cite{alex}, \cite{alex2}), one obtains a version
of the case $\dim X = \dim B + 2$. Here each fiber is a
semi-log-canonical surface.

Up until recently, not much has been known about the case $\dim X
>\dim B+2$.  Often one finds remarks of the following flavor: ``since
we do not have a semistable reduction result over a base of higher
dimension, we will work around it in the following technical
manner...''.

\subsection{Definition of semistable families} We  give here a
description of the best possible kind of morphisms we have in mind.

Let $f:X\to B$ be a flat morphism of nonsingular projective varieties
with connected fibers.  Somewhat informally, we say that $f$ is {\bf
semistable} if for each point $x\in X$ with $f(x) = b$ there is a
choice of formal coordinates $\hat{B}_b = \spec\ k[[t_i]]$ and $\hat{X}_x =
\spec\
k[[x_j]]$, such that $f$ is given by: $$t_i =
\prod_{j=l_{i-1}+1}^{l_i} x_j. $$ Here $0 = l_0 < l_1 \cdots < l_m
\leq n$, where $n=\dim X$ and $m = \dim B$.  To be more precise, we
give things a more global structure using the notion of a toroidal
morphism. At the same time we describe a slightly weaker condition
which will appear below:

\begin{dfn}
The morphism $f:X\to B$ above is called {\bf weakly semistable} if
\begin{enumerate}
\item the varieties $X$ and $B$ admit toroidal structures $U_X\subset
X$ and $U_B\subset B$, with $U_X=f^{-1}U_B$;
\item with this structure, the morphism $f$ is toroidal;
\item the morphism $f$ is equidimensional;
\item all the fibers of the morphism $f$ are reduced;
and
\item $B$ is nonsingular.
\end{enumerate}
If also $X$ is nonsingular, we say that the morphism $f:X\to B$ is
 {\bf semistable}.
\end{dfn}

\subsection{The ultimate goal} The result one would really like to have is:

\begin{conj}\label{conj-semistable}
Let $X\rar B$ be a surjective morphism of complex projective
varieties with geometrically integral generic fiber. There is a projective
alteration $B_1\rar B$, and a projective
modification $Y\rar X\times_B B_1$, such that $Y\rar B_1$ is semistable.
\end{conj}

Na\"{\i}vely one might hope to have each fiber isomorphic to a divisor
of normal crossings. But already in the case of a 2-parameter family
of surfaces $t_1 = x_1x_2;\, t_2 = x_3 x_4$, this is impossible. It
seems that the definition above is the best one can hope for.

\subsection{A. J. de Jong's results} In \cite{dj}, Johan de Jong shows,
among many other results, that if one allows $Y\rar X\times_B B_1$ to
be an alteration instead of a modification, one can make $Y\rar B_1$
very nice indeed: $Y$ is nonsingular, $Y\rar B_1$ is semistable as in
the definition above, and moreover it can be written as a composition
of nodal curve fibrations $Y=Y_0\to Y_1\to \cdots \to Y_k=B_1$.

De Jong's methods and ideas will serve as a starting point for
investigating the semistable reduction conjecture.

\subsection{Our main result} The main result of this paper is the
following:

\begin{th}[Weak semistable reduction]\label{th-weak-semistable-reduction}
Let $X\rar B$ be a surjective morphism of complex projective
varieties with geometrically integral generic fiber. There exist an
alteration $B_1\rar B$ and a modification
$Y\rar X\times_B B_1$, such that $Y\rar B_1$ is weakly semistable.
\end{th}

 With a little more work we will get $X$ to have only quotient
singularities. There are many cases (such as when $f$ is a family of
surfaces) where we can actually prove the semistable reduction
conjecture. These will be pursued elsewhere. Hopefully, by the time
this paper achieves its final form the conjecture will be fully
proven.

\subsection{Mild morphisms} A few words are in order about the significance of
our result. Note that the property of a morphism being semistable is
far from being stable under base changes. One may ask, what remains
from semistability after at least {\em dominant} base changes? Here is
a suggestion:

\begin{dfn}
We define a morphism $X\rar B$ as above to be {\bf mild}, if for any
dominant $B_1\rar B$ where $B_1$ has at most rational Gorenstein
singularities, we have that $X\times_B B_1$ has at most rational
Gorenstein singularities as well.
\end{dfn}

Mild morphisms arise naturally in moduli theory. Indeed, mild families
of curves are precisely nodal families; families of Gorenstein
semi-log-canonical surfaces mentioned above are mild. For a discussion
of why mild morphisms are useful, see \cite{fibered}. In fact, the
paper \cite{fibered} would have been much simplified, had mild
reduction been available.

Already in the case $\dim B=1$, mild reduction is a much easier task
than semistable reduction. Indeed, lemma 2 on page 103 of \cite{te},
and the discussion there, already give mild reduction in this
case. The delicate combinatorics of chapter III of \cite{te} is not
used for this purpose.

It will be shown (see section \ref{mildness}) that weakly semistable morphisms
are indeed mild.

\subsection{Structure of the proof}

After the introduction, section \ref{toroidal-morphisms} will be
devoted to a general discussion of toroidal morphisms. The proof
itself will begin with section \ref{toroidal-reduction}.

Semistable reduction has at least two flavors: first, the fibers of
the morphism $Y\rar B_1$ should have nice local defining
equations. Second, the family should have nice algebraic
properties. We will perform a number of reduction steps, incrementally
improving one or the other of these flavors.

\subsubsection{Toroidal reduction} In the first step,
carried out in section \ref{toroidal-reduction}, we will show that any
morphism can be modified to a toroidal morphism. The construction is
inspired by the inductive procedure of \cite{dj}, and follows closely
the proofs in \cite{aj}.

Just as in \cite{aj}, the construction we give is very
non-canonical. Even when the generic fiber of $X\to B$ is smooth, it
will be blown up during the construction. One hopes that methods such
as those of \cite{bm} or \cite{villa} could be adapted to this
situation and give a more canonical procedure.

It is tempting to state the following conjecture.

\begin{conj} Let $X\to B$ be a morphism as in the theorem. Let $U\subset B$ be
an open set over which $X$ is toroidal, and let $\Sigma = B \setmin U$. There
exists modifications $X'\to X$ and $B'\to B$, each of which is the composition
of a sequence of blowings up with smooth centers lying over $\Sigma$,  and a
lifting $X'\to B'$ which is toroidal.
\end{conj}

It should be noted, that in view of recent results of Morelli \cite{m} and
W{\l}odarczyk \cite{w}, this conjecture implies the strong blow-up - blow-down
conjecture.

\subsubsection{Improving the toroidal morphism} In sections
\ref{remove-horizontal} and \ref{equidimensional} we perform a couple
of simple reduction steps to improve our situation. {Let $f:(U_X
\subset X) \rar(U_B\subset B)$ be any toroidal morphism, with $B$
nonsingular. By the results of \cite{te}, we can find a toroidal
resolution of singularities $X'\rar X$. Let $f':X'\rar B$ be the
resulting projection. We first show that now $f^{-1} U_B\subset X'$ is
also a toroidal embedding, which is easier to handle: there are no
horizontal divisors. For convenience, we replace $X\rar B$ by the new
morphism.  We remark that one can proceed a fair distance without
removing these horizontal divisors, and, we believe, the results one
can obtain are of interest (e.g., the inductive structure of de Jong
can be preserved), but this would make the present paper much more
cumbersome, so we delay that investigation to a future occasion.

Now, our morphism $X\rar B$ is not necessarily equidimensional. We
repair this by an appropriate decomposition of the associated conical
polyhedral complexes $\Delta_X$ and $\Delta_B$. We make sure that,
after the modification, the base remains nonsingular, and then the morphism
is automatically flat.

\subsubsection{Kawamata's trick and reduced fibers.}  We start  section
\ref{reduced-fibers} with a discussion of Kawamata's trick and its
relation with toroidal morphisms in some detail. Then we use
Kawamata's trick to find a finite base change, after which all the
fibers are reduced. This finishes the proof of the main theorem, since
the resulting morphism is weakly semistable. A variant of Kawamata's trick for
global ``index 1 covers'' is discussed in section \ref{sec-cartier}.

\subsubsection{Mild reduction.}
We begin section \ref{mildness} by checking that the resulting fibers
are Gorenstein.

 Using a base change and descent argument, and the fact that toroidal
singularities are always rational, we then prove that the resulting
family is mild.

\subsubsection{Combinatorial restatement} In section \ref{combinatorial} the
semistable reduction conjecture is restated purely in combinatorial terms. We
end the paper with a discussion of the problems one encounters when
trying to go from weak semistable reduction to semistable reduction.

\subsection{Acknowledgments} We would like to thank O. Gabber, A.J. de Jong,
H. King and K. Matsuki for helpful and inspiring discussions.

\subsection{Terminology}

A {\bf modification} is a proper birational morphism of irreducible
varieties.

An {\bf alteration} $a:B_1\rar B$ is a proper, surjective, generically
finite morphism of irreducible varieties, see \cite[2.20]{dj}.  The
alteration $a$ is a {\bf Galois alteration} if there is a finite group
$G\subset \Aut_B (B_1)$ such that the associated morphism $B_1/G\rar
B$ is birational, compare \cite[5.3]{dj2}.

\section{Toroidal morphisms}\label{toroidal-morphisms}

We have collected in this section some notations and preliminaries
about toric varieties, toroidal embeddings, and their morphisms (see
\cite{te} for details).

\begin{rem}
Our approach here is based on the formalism of \cite{te}. A different approach,
using {\em logarithmic structures}, was developed by K. Kato, see \cite{kato},
\cite{kato1}. It is our belief, that the approach via logarithmic structures
should eventually prevail - it provides us with a flexible category, in which
toroidal embeddings ($=$ logarithmically regular schemes) and toroidal
morphisms ($=$ logarithmically smooth(!) morphisms) play a special role.
Some of our statements below are rendered almost trivial with Kato's formalism,
e.g. Lemmas \ref{lem-tor-composition} and \ref{lem-tor-prod}.

The reason  we decided to stick with the formalism of \cite{te} is, that the
theory of logarithmic structures is not yet in stable form (see the many
flavors of such structures introduced in Kato's papers), and, more importantly,
it has not yet gained widespread acceptance as a basic formalism. It might have
turned away some readers (especially those combinatorially inclined) had we
used the theory of logarithmic structures throughout.

It is also worth noting, that Kato's notion of a fan, although it has a nice
structural morphism, is much less amenable to combinatorial manipulation than
the polyhedral complexes of \cite{te}.
\end{rem}

\subsection{Toric varieties} Given a lattice $N\cong\bfz^n$, its dual
$M=\Hom(N,\bfz)$,
a strictly convex rational polyhedral cone
$\sigma\subset N_\bfr = N\otimes\bfr$ with its dual $\sigma^\vee =
\{m\in M_\bfr | m(u)\geq 0$ for all $u\in \sigma\}$, we define the {\bf affine
toric variety} $X_\sigma = \spec S[\sigma]$ where
$S[\sigma]$ is the semigroup algebra of $\sigma^\vee\cap M$ over the
ground field. If more than one toric variety is considered, we use a
subscript: $N_\sigma$, $M_\sigma$.

We denote by $\sigma^{(1)}$ the 1-dimensional edges of $\sigma$.
The indivisible points $v$ in $\sigma^{(1)}\cap N$ are called the {\bf
primitive
points} of $\sigma$. The variety $X_\sigma$ is nonsingular if and only if the
primitive points of $\sigma$ form a part of a basis of $N$. In that
case we say that $\sigma$ is nonsingular.

The toric variety $X_\sigma$ contains an n-dimensional algebraic torus
$T=\bfg^n_m$ as an open dense subset, and the action of $T$ on itself
extends to an action on $X_\sigma$. Thus, $X_\sigma$ is a disjoint
union of orbits of this action. There is a one-to-one correspondence
between the orbits and the faces of $\sigma$. In particular,
1-dimensional faces $\bfr_{+}v_i$ correspond to codimension~1 orbits
$\bfo_{v_i}$.

A {\bf toric morphism} $f:X_\sigma\rar X_\tau$ is a dominant equivariant
morphism of toric varieties defined by a linear map $f_\Delta:
(N_\sigma,\sigma) \rar (N_\tau,\tau)$. We use the same notation for
the scalar extension $f_\Delta: N_\sigma\otimes\bfr \rar
N_\tau\otimes\bfr$.

\subsection{Toroidal embeddings}
Given a normal variety $X$ and an open subset $U_X\subset X$,
the embedding $U_X\subset X$ is called {\bf toroidal} if for every
closed point $x\in X$ there exist a toric variety $X_\sigma$, a point
$s\in X_\sigma$, and an isomorphism of complete local algebras
\[ \hat{\co}_{X,x} \cong \hat{\co}_{X_{\sigma},s} \]
so that the ideal of $X\setmin U_X$ corresponds to the ideal of
$X_{\sigma}\setmin T$. Such a pair $(X_\sigma, s)$ is called a local
model at $x\in X$. By restricting $X_{\sigma}$ if necessary, we can
assume that the orbit of $s$ is the unique closed orbit in
$X_{\sigma}$.

\begin{dfn}\label{def-toroidal-map}
A dominant morphism ${f}:(U_X\subset X)\rar(U_B\subset B)$ of toroidal
embeddings is called {\bf toroidal} if for every closed point $x\in X$
there exist local models $(X_{\sigma},s)$ at $x$, $(X_{\tau},t)$ at
${f}(x)$ and a toric morphism $g: X_{\sigma}\rar X_{\tau}$ so that the
following diagram commutes
\[
\begin{CD}
\hat{\co}_{X,x} @<{\cong}<< \hat{\co}_{X_{\sigma},s}\\
@A{\hat{f}^{*}}AA @AA{\hat{g}^{*}}A\\ \hat{\co}_{B,f(x)} @<{\cong}<<
\hat{\co}_{X_{\tau},t}
\end{CD}
\]
where $\hat{f}^{*}$ and $\hat{g}^{*}$ are the ring homomorphisms
induced by $f$ and $g$.
\end{dfn}

\subsection{Cones and polyhedral complexes} Let $X\setmin U_X = \cup_{i\in
I}E_i$ where $E_i$ are irreducible
and have codimension~1. We will assume that all the $E_i$ are normal,
that is, $U_X\subset X$ is a {\bf toroidal embedding without
self-intersection} (also known as a {\bf strict} toroidal embedding). In that
case, we can use the irreducible components of $\cap_{i\in J}
E_i$ for all $J\subset I$ to define a stratification of $X$ (these components
are the closures of strata). Closures of strata
formally corresponds to closures of orbits in local models. Since a toric
morphism maps orbits to orbits, a toroidal morphism maps strata to
strata.

Let Y be a stratum in $X$, which is by definition an open set in an irreducible
component of $\cap_{i\in J} E_i$ for some $J\subset I$. The star of
$Y$ is the union of strata in whose closure $Y$ lies (each of these corresponds
to some $K\subset  J\subset I$). To
the stratum $Y$ we associate
\begin{enumerate}
\item $M^Y$: --  the group of Cartier divisors in $\Star(Y)$ supported
in $\Star(Y)\setmin U_X$
\item $N^Y$: -- $\Hom(M^Y,\bfz)$
\item $M_{+}^Y \subset M^Y$: -- effective Cartier divisors
\item $\sigma^Y\subset N^Y_\bfr$: -- the dual of $M_{+}^Y$
\end{enumerate}
It is shown in \cite{te} (Corollary 1, page 61) that if
$(X_{\sigma},s)$ is a local model at $x\in X$ in the stratum $Y$, then
\begin{eqnarray*}
M^Y &\cong& M_\sigma/\sigma^{\bot} \\ \sigma^Y &\cong& \sigma
\end{eqnarray*}
The cones $\sigma^Y$ glue together to form a polyhedral complex
$\Delta_X = (|\Delta_X|,\{\sigma^Y\},\{M^Y\})$, where the lattices $M^Y$ form
an
integral
structure on $\Delta_X$. Equivalently, instead of $M^Y$ we may give
the lattices $N^Y$ and embeddings $\sigma^Y\hookrightarrow
N^Y_\bfr$. We also denote $M_{\sigma^Y}=M^Y$,
$N_{\sigma^Y}=N^Y$. Then, comparing to the lattices of local models,
\begin{eqnarray*}
M_{\sigma^Y} &\cong& M_\sigma/\sigma^{\bot} \\ N_{\sigma^Y} &\cong&
N_\sigma\cap \mbox{Span}(\sigma) \\
\end{eqnarray*}
If $N_{\sigma^Y}\neq N_\sigma$ we can choose a splitting of $N_\sigma$
so that at a point $x\in Y$ the local model is in the form
$(X_{\sigma'}\times\bfg_m^l,(s',1))$ where $(N_{\sigma'},\sigma')\cong
(N_{\sigma^Y},\sigma^Y)$.

\begin{lem}
A toroidal morphism ${f}:X\rar B$ induces a morphism
$f_\Delta:\Delta_X \rar
\Delta_B$, which for each cone $\sigma^Y$ is the
restriction of $g_\Delta: (\sigma,N_\sigma)\rar(\tau,N_\tau)$ where
$(X_{\sigma},s)$, $(X_{\tau},t)$ are local models at $x\in Y\subset
X$, $f(x)\in Z\subset B$, and $g_\Delta$ is the linear map determined
by the toric morphism $g: X_{\sigma}\rar X_{\tau}$ in
Definition~\ref{def-toroidal-map}.
\end{lem}

{\bf Proof.} To see that the maps $f_\Delta$ defined for different
cones $\sigma^Y$ agree on the overlaps it suffices to notice that the
dual morphism $f_\Delta^\vee:M_{\tau^Z}\rar M_{\sigma^Y}$ is defined
by pulling back a Cartier divisor and restricting it to
$\Star(Y)\setmin U_X$. Since the pullback is defined independently of
the stratum, we see that $f_\Delta$ is well defined. \qed

\begin{rem}
Note that the polyhedral morphism $f_\Delta:\Delta_X \rar
\Delta_B$ is well defined even if $f$ is not toroidal, as long as
$f(U_X)\subset U_B$.
\end{rem}

\begin{lem}\label{lem-tor-composition}
If $e:X\rar Y$ and $f:Y\rar Z$ are toroidal morphisms, then $f\circ
e:X\rar Z$ is also toroidal.
\end{lem}

\begin{rem} This lemma is a triviality if one uses logarithmic structures.
\end{rem}

{\bf Proof.}  Let $x\in X$, $y=e(x)\in Y$, $z=f(y)\in Z$, and choose
local models at $x$, $y$ and $z$ as in
Definition~\ref{def-toroidal-map}. Consider the tower
\[
\begin{CD}
\hat{\co}_{X,x} @<{\cong}<< \hat{\co}_{X_{\sigma},s}\\
@A{\hat{e}^{*}}AA @AA{\hat{g}^{*}}A\\ \hat{\co}_{Y,y} @<{\cong}<<
\hat{\co}_{X_{\tau_1},t_1}\\ @A{id}AA @AA{\alpha}A\\ \hat{\co}_{Y,y}
@<{\cong}<< \hat{\co}_{X_{\tau_2},t_2}\\ @A{\hat{f}^{*}}AA
@AA{\hat{h}^{*}}A\\ \hat{\co}_{Z,z} @<{\cong}<< \hat{\co}_{X_{\rho},r}
\end{CD}
\]
where the upper and lower squares commute by the definition of
toroidal morphism, and where $\alpha$ is defined by tracing the other
three sides of the middle square. Then the middle square also
commutes.

Since $\tau_1$ and $\tau_2$ are isomorphic, we can take
$\tau_1=\tau_2=\tau$ and $X_{\tau_1}= X_{\tau_2}= X_{\tau}$. The map
$\alpha$, of course, need not be the identity. Let the coordinate
rings of the tori in $X_{{\sigma}}$, $X_{{\tau}}$ and $X_{{\rho}}$ be
$k[x_1,x_1^{-1},\ldots,x_n,x_n^{-1}]$,
$k[y_1,y_1^{-1},\ldots,y_m,y_m^{-1}]$ and
$k[z_1,z_1^{-1},\ldots,z_l,z_l^{-1}]$, respectively, so that the toric
morphisms $g$ and $h$ are defined by
\[ y_i \mapsto \prod_{j=1}^{n} x_j^{a_{ij}}, \hspace{.5in}
   z_i \mapsto \prod_{j=1}^{m} y_j^{b_{ij}} \] for some $a_{ij},
b_{ij} \in \bfz$.

The maps in the third line of the tower identify the group of
T-invariant Cartier divisors in $X_{{\tau_2}}$ with the group of
Cartier divisors in $Y$ supported in $Y\setmin U_Y$ and passing
through $y$. The maps in the second line of the tower identify the
latter group with the group of T-invariant Cartier divisors in
$X_{{\tau_1}}$. Hence $\alpha$ induces a group homomorphism between
T-invariant Cartier divisors of $X_{{\tau_2}}$ and
$X_{{\tau_1}}$. Since T-invariant Cartier divisors in $X_{{\tau}}$ are
given by products of $y_i$, then (using the same letter $y_i$ for the
image of $y_i$ in the completed local ring) $\alpha$ maps
\[ y_i \mapsto u_i \prod_{j=1}^{m} y_j^{c_{ij}} \]
where $u_i$ are units in $\hat{\co}_{X_{{\tau}},t_1}$, and $c_{ij}\in
\bfz$.

The composition $\hat{g}^{*}\circ\alpha\circ\hat{h}^{*}$ is then
defined by
\[ z_i \mapsto v_i \prod_{j=1}^{m} x_i^{d_{ij}} \]
where the matrix with entries $d_{ij}$ is the product of the matrices
with entries $a_{ij}$, $c_{ij}$ and $b_{ij}$, and where
$v_i=\hat{g}^{*}(\prod_j u^{c_{ij}})$ are units in
$\hat{\co}_{X_{{\sigma}},s}$. The matrix $(d_{ij})$ is equivalent to a
matrix in the form
\[
\begin{pmatrix} D \\ 0\end{pmatrix}
\]
where D is an $l\times l$ diagonal matrix with diagonal entries
$d_1,\ldots,d_l \in \bfz$. Hence we can change the coordinate
functions $x_i$ and $z_i$ of the respective tori so that
$\hat{g}^{*}\circ\alpha\circ\hat{h}^{*}$ maps
\[  z_i \mapsto v_i x_i^{d_i}\]
If we now set
\[ \tilde{x_i} = v_i^{\frac{1}{d_i}} x_i \]
then in these new coordinates the composition
$\hat{g}^{*}\circ\alpha\circ\hat{h}^{*}$ is induced by a toric
morphism defined by
\[ z_i \mapsto x_i^{d_i} \]
and thus $f\circ e:X\rar Z$ is toroidal.  \qed

Given a toroidal embedding $U_X\subset X$ with polyhedral complex
$\Delta_X$, and a subdivision $\Delta'_X$ of $\Delta_X$, one
constructs (see \cite{te}) a new toroidal embedding $U_{X'}\subset X'$ with
polyhedral complex $\Delta'_X$, and a birational toroidal morphism $f': X'\rar
X$ such that the induced map of the polyhedral complexes $ \Delta'_X \rar
\Delta_X$ is the given subdivision. If $Y$ is a stratum in $X$ corresponding
to the cone $\sigma^Y \in \Delta_X$, and if $\sigma'\subset \sigma$ is
a cone in the subdivision, define
\[ V_{\sigma'} = \spec_{\Star(Y)} \sum_{D\in {\sigma'}^{\vee} \cap M^Y}
\cO_{\Star(Y)}(-D), \]
where the sum is taken inside the field of rational functions of
$\Star(Y)$.
Then $X'$ is formed by gluing together the open sets $V_{\sigma'}$.

A subdivision $\Delta_X'$ of $\Delta_X$ is called {\bf projective} if
there exists a continuous  function  $\psi:|\Delta_X|\rar \bfr$, taking
rational values on $\sigma\cap N_\sigma$,  which is
convex and piecewise linear on each cone $\sigma\in\Delta_X$, and the
largest pieces in $\sigma$ where $\psi$ is linear are the cones of the
subdivision. Such $\psi$ is called a {\bf good function} (or {\em lifting
function}, or {\em order function}), for the
subdivision $\Delta_X'$ of $\Delta_X$. A projective subdivision
corresponds to a projective modification $X'\rar X$.

\begin{lem}\label{lem-lifting} Let $f:X\to B$ be a toroidal morphism,
$f_\Delta:\Delta_X\to
\Delta_B$ the associated morphism of polyhedral complexes. Let $X'\to
X$ and $B'\to B$ be toroidal modifications, with associated
subdivisions $\Delta_{X'}$ and $\Delta_{B'}$. Then $f$ lifts to a
morphism $f': X' \to B'$ if and only if for each cone $\sigma' \in
\Delta_{X'}$, there exists a cone $\tau' \in \Delta_{B'}$ such that
$f_\Delta(\sigma') \subset \tau'$.
\end{lem}

{\bf Proof.} Let $\sigma'\subset \sigma^Y \in \Delta_X$ and
$\tau'\subset\tau^Z\in\Delta_B$ be cones in the subdivisions such that
$f_\Delta(\sigma') \subset \tau'$. The homomorphism
$\co_{\Star(Z)} \rar \co_{\Star(Y)}$ extends to
\[ \sum_{D\in {\tau'}^{\vee} \cap M^Z} \cO_{\Star(Z)}(-D) \rar
\sum_{E\in {\sigma'}^{\vee} \cap M^Y} \cO_{\Star(Y)}(-E) \]
because $f^{\vee}_{\Delta} ({\tau'}^\vee) \subset
{\sigma'}^{\vee}$. This shows that the rational map $f'$ is a morphism
on each $V_{\sigma'}$. Since $V_{\sigma'}$ cover $X'$, $f'$ is a
morphism.

Conversely, suppose $f'$ is a morphism. Let $x'\in X', f'(x')=b'\in
B'$. Since the three morphisms $X'\rar X$, $X\rar B$ and $B'\rar B$
are toroidal, we have a toric morphism $g: X_{\sigma'}\rar X_{\tau'}$
where $(X_{\sigma'},s')$ and $(X_{\tau'},t')$ are local models at $x'$
and $b'$. Thus $g_{\Delta}(\sigma') \subset \tau'$. But
$g_{\Delta}|_{\sigma'} = f_{\Delta}|_{\sigma'}$, hence
$f_{\Delta}(\sigma') \subset \tau'$. Since this is true for $x'$ in
any stratum of $X$, we get that for any cone $\sigma' \in
\Delta_{X'}$, $f_\Delta(\sigma') \subset \tau'$ for some $\tau' \in
\Delta_{B'}$. \qed

\section{Toroidal reduction}\label{toroidal-reduction}
\subsection{Statement of result}

The purpose of this section is to modify any family of varieties into
a toroidal morphism.

\begin{th}\label{th-toroidal-reduction}
Let $f:X\to B$ be a projective, surjective morphism with geometrically
integral generic fiber, and assume $B$ integral. Let $Z\subset X$ be a proper
closed subscheme. There exist a diagram as follows:
$$\begin{array}{lclcl} U_X & \subset & X' &\stackrel{m_X}{\to} &X \\
                     \dar & & \dar f' & & \dar f\\ U_B & \subset & B'
                     & \stackrel{m_B}{\to} &B \end{array}
$$
such that $m_B$ and $m_X$ are modifications, $X'$ and $B'$ are
nonsingular, the inclusions on the left are toroidal embeddings, and
such that
\begin{enumerate}
\item $f'$ is toroidal.
\item Let $Z' = m_X^{-1}Z$. Then $Z'$ is a strict normal crossings
divisor, and $Z'\subset X'\setmin U_{X'}$.
\end{enumerate}
\end{th}

\subsection{To begin the proof,} we proceed by induction on the relative
dimension of $f$.

If the relative dimension of $f$ is 0, let $m_X:X'\to X$ be a
resolution of singularities such that $Z'=m_X^{-1}Z$ is a strict
normal crossings divisor, let $B'=X'$ and $m_B = f\circ m_X$, and
$f'=id$ the identity.

Assume we have proven the result for morphisms of relative dimension
 $n-1$, and consider the case $\reld f = n$.

\subsection{Preliminary reduction steps}
First, we may replace $X$ by its normalization, therefore we may
assume $X$ normal, and by blowing up $Z$ in $X$ we may assume $Z$ a
Cartier divisor. Let $\eta\in B$ be the generic point of $B$. By the
projectivity assumption we have $X\subset \bfp^N_B$ for some
$N$. Choosing a generic projection $\bfp^N_\eta \das \bfp^{n-1}_\eta$
we get a rational map $X_\eta \das \bfp^{n-1}_\eta$. Replacing $X$ by
the closure of the graph of this map, we may assume that we have a
morphism $g:X \to \bfp^{n-1}_B = P$.

\subsection{Semistable reduction of a family of curves}
By \cite{dj2}, Theorem 2.4, we have a diagram as follows:

$$\begin{array}{lcl} X_1 & \stackrel{\alpha}{\rar} & X \\ \dar g_1 & &
                     \dar g \\ P_1 & \stackrel{a}{\rar} & P \\ & &
                     \dar \\ & & B \end{array}
$$
 and a finite group $G\subset \Aut_PP_1$, with the following
 properties:

\begin{enumerate}
\item The morphism $a:P_1\to P$ is a Galois alteration with Galois
group $G$.
\item The action of $G$ lifts to $\Aut_XX_1$, and $\alpha:X_1\to X$ is
a Galois alteration with Galois group $G$.
\item There are $n$ disjoint sections $\sigma_i:P_1\rar X_1$ such that
the strict altered transform $Z_1\subset X_1$ of $Z$ is the union of
their images, and $G$ permutes the sections $\sigma_i$.
\item The morphism $g_1: X_1 \to P_1$ is a nodal family of curves, and
$\sigma_i(P_1)$ is disjoint from $\operatorname{Sing}g_1$.
\end{enumerate}

We may replace $X$, $P$ and $Z$ by $X_1/G$ and $P_1/G$, and
$\alpha^{-1}Z/G$.  Note that $\alpha^{-1}Z/G$ is not necessarily equal
to the union of the images of $\sigma_i$, but the complement lies over
a proper closed subset in $P_1$.

\subsection{Using the inductive hypothesis}
Let $\Delta\subset P$ be the union of the loci over which $Z,P_1$ or
$X_1$ are not smooth. We apply the inductive assumption to
$\Delta\subset P \to B$, and obtain a diagram as follows:

$$\begin{array}{lclcl} U_P & \hookrightarrow & P' &\stackrel{m}{\to} &P \\
\dar                      &                      & \dar & & \dar\\
 U_B & \hookrightarrow & B' &\to &B
                     \end{array}
$$

Such that $P', B'$ are nonsingular, $P' \to P$ and $B' \to B$ are
modifications, the left square is a toroidal morphism, and
$m^{-1}\Delta$ is a divisor of strict normal crossings contained in
$P' \setmin U_P$.

We may again replace $P, B$ by $P', B'$, and further we may replace
$X,X_1, P_1,Z$ and $ \sigma_i$ by their pullback to $P'$. In
particular $P\to B$ has a toroidal structure, and $P_1\to P$ is
unramified over $U_P$. By Abhyankar's lemma, since $P_1$ is normal, it inherits
a toroidal structure given by $U_{P_1} = m^{-1}U_P$ as well, so that $P_1\to P$
is a toroidal finite morphism.

\subsection{Conclusion of proof} Now $X_1\to P_1$ is a nodal family which is
smooth over $U_P$, therefore it as well inherits a toroidal structure
$U_{X_1}\subset U_X$, where $U_{X_1} = ({g_1}^{-1} U_{P_1}) \setmin
(\cup\sigma_i(P_1))$; e.g. local equations around a node are of the form
$ uv = f(t)$, where $f(t)$ is a monomial on $P_1$. Notice that
$\alpha^{-1} Z $ is a divisor contained in $U_{X_1}$ (see
\cite{aj}, 1.3).

In this situation we can apply the procedure of \cite{aj}, section 1.4
to make the group $G$ act toroidally on $X_1$: first we blow up the
scheme $\operatorname{Sing} g_1$ to separate the branches of the
nodes. Then we are in the situation of Proposition 1.8 of \cite{aj},
namely there is a canonical $G$-equivariant blowup $d:\tilde{X_1} \to
X_1$ such that $G$ acts strictly toroidally on
$b^{-1}U_{X_1}\subset\tilde{X_1}$. Let ${X'} =\tilde{X_1}/G$, then $
{X'}\to B $ inherits a toroidal structure and ${X'}\to X$ is
birational; moreover, ${Z'}\subset {X'}$ is a divisor contained in
${X'}\setmin U_{{X'}}$. Applying toroidal resolution of singularities,
the induction step is proven. \qed

\section{Removing horizontal divisors}\label{remove-horizontal}

We may now replace $X\to B$ by $X'\to B'$, and thus we may assume that the
morphism $f$ is toroidal.
Our goal in this section is to arrive at a situation where $f^{-1} U_B = U_X$.

The rough idea is, that a morphism between nonsingular toroidal
embeddings $f:X\to B$ is locally given by monomials $t_i = x_1^{k_1}
\cdots x_r^{k_r}$, in which the variables defining horizontal divisors
cannot appear, so these divisors are unnecessary in the toroidal
description. We make this precise by a simple translation argument.

\begin{prp}
Let $U_X\subset X$ and $U_B \subset B$ be nonsingular toroidal
embeddings and $f: X\rar B$ a surjective toroidal morphism. Then,
denoting $U_X'= f^{-1}(U_B) \supset U_X$, we have that $U_X'\subset X$ is a
toroidal embedding, and $f: (U_X'\subset X)\rar (U_B\subset B)$ is a
toroidal morphism.
\end{prp}

{\bf Proof.}  Since $f$ maps $U_X$ into $U_B$, $f^{-1}(B\setmin
U_B)$ as a set is a union of divisors supported in $X\setmin
U_X$. In local models these divisors are all T-invariant.

Consider local models $(X_{\sigma},s)$ at $x$, $(X_{\tau},t)$ at
$f(x)$, and the toric morphism $g:X_\sigma\rar X_\tau$. We may assume
that $v_1,\ldots,v_n$ is a basis of $N_\sigma$ and $\sigma$ is
generated by $v_1,\ldots,v_k$. Then $X_\sigma \cong
\bfa^k\times\bfg_m^{n-k}$, and we may take $s=(0,1)$. Let the closures
of the orbits corresponding to $v_1,\ldots,v_j$ be the horizontal
divisors. That means, $g_\Delta(v_i)=0$ for $i=1,\ldots,j$, and $g$
factors through the projection:
\[ g: X_\sigma \cong \bfa^j\times\bfa^{k-j}\times\bfg^{n-k}_m \rar
\bfa^{k-j}\times\bfg^{n-k}_m \to X_\tau \]
Now take $s'=(1,0,1)\subset
\bfa^j\times\bfa^{k-j}\times\bfa^{n-k}$. From the factorization of $g$
we see that $g(s')=g(s)$. Translation by $(1,0,0)$ gives an
isomorphism of the local rings at $s$ and $s'$ so that the ideals of
the T-invariant divisors passing trough $s'$ corresponds to the ideal
of the vertical T-invariant divisors passing through $s$. Thus,
$(X_\sigma,s')$ is a local model for the embedding $U_X'\subset X$ at
$x\in X$ and $g:(X_\sigma,s')\rar(X_\tau,t)$ is the toric morphism of
the local models representing $f$. \qed

\section{Making the morphism equidimensional}\label{equidimensional} The goal
of this section is to perform modifications on $B$ and $X$, after which
the morphism becomes equidimensional. First a lemma which
characterizes equidimensional toroidal morphisms:

\begin{lem} Let $f:X\to B$ be a surjective toroidal morphism, $f_\Delta:
\Delta_X\to\Delta_B$ the associated morphism of polyhedral complexes.
Then $f$ is equidimensional if and only if for any cone $\sigma \in
\Delta_X$, we have $f_\Delta(\sigma) \in \Delta_B$. That is, the image
of a cone of $\Delta_X$ is a cone of $\Delta_B$.
\end{lem}

{\bf Proof.} Computing the dimension of a local ring commutes with
taking the completion. Thus, it suffices to consider local models. The
generic fiber of a toric morphism $f: X_\sigma \rar X_\tau$ has
dimension $\rk(N_\sigma) - \rk(N_\tau)$. Now $f$ maps a
$k$-dimensional orbit corresponding to a
$(\rk(N_\sigma)-k)$-dimensional face $\sigma'$ of $\sigma$ onto some
$l$-dimensional orbit corresponding to a $(\rk(N_\tau)-l)$-dimensional
face $\tau'$ of $\tau$. Thus $f$ is equidimensional if and only if
\[ (\rk(N_\sigma)-k) - (\rk(N_\tau)-l) \leq \rk(N_\sigma) - \rk(N_\tau), \]
that means $l \leq k$. Every $k$-dimensional cone maps to an
$l$-dimensional cone for $l\leq k$ if and only if the image of every
cone is a cone. \qed

\begin{rem} In case $\tau$ is simplicial the condition of the theorem
is equivalent to the statement that all 1-dimensional faces of
$\sigma$ map to 0 or 1-dimensional faces of $\tau$.
\end{rem}

The following lemma is a slight generalization of the toric Chow's
lemma (\cite{danilov} 6.9.2 page 119).

\begin{lem} Given a polyhedral complex $\Delta$ and a subdivision
$\Delta'$ of $\Delta$, there exists a projective subdivision
$\Delta''$ of $\Delta$ which refines $\Delta'$.
\end{lem}

{\bf Proof.}  First, we show that it suffices to find a ``good''
function $\psi$ (see Section~\ref{toroidal-morphisms}) on each cone
$\sigma\in\Delta$. Indeed, if for each $\sigma\in\Delta$ we have found
a good function $\psi_\sigma: \sigma\rar\bfr$, then a good function on
$|\Delta|$ is the sum
\[ \psi = \sum_{\sigma\in\Delta} \bar{\psi}_\sigma \]
where $\bar{\psi}_\sigma: |\Delta|\rar\bfr$ is a good extension of
$\psi_\sigma$ constructed as follows. To extend $\psi_\sigma$ to a
cone $\tau\in\Delta$ we proceed by induction on the dimension of
$\tau$ ({\em cf.} \cite{te} Lemma~1, page~33). If $\dim \tau =
1$ and $\tau$ is not a face of $\sigma$, define
$\bar{\psi}_\sigma|_\tau\equiv 0$. Now assume that $\dim {\tau} >
1$ and $\bar{\psi}_\sigma$ is defined on $\partial\tau$. Choose a
point $x$ in the relative interior of $\tau$ and define
\[ \bar{\psi}_\sigma (\lambda x + (1-\lambda) y) = \lambda C +
(1-\lambda) \bar{\psi}_\sigma(y), \hspace{0.3in} y\in\partial\tau,
0\leq\lambda\leq 1.\]
For big $C\in\bfq$ the extension $\bar{\psi}_\sigma|_\tau$ is convex.

Now let $\sigma\in\Delta$ with $\dim {\sigma}=n$. For every $n-1$
dimensional cone $\tau\in\Delta'$ in the subdivision of $\sigma$ choose
a linear rational function $l_\tau: \sigma\rar\bfr$ such that $\tau$
is the subset of $\sigma$ defined by $l_\tau=0$. Then the sum
\[ \psi_\sigma(x) = -\sum_{\tau}|l_\tau(x)| \]
is a good function on $\sigma$. \qed

\begin{prp}
Let $U_X\subset X$ and $U_B \subset B$ be toroidal embeddings with
polyhedral complexes $\Delta_X$ and $\Delta_B$ respectively, and
assume that $B$ is nonsingular. Let $ f: X\rar B$ be a surjective
toroidal morphism. Then there exist projective subdivisions
$\Delta_X'$ of $\Delta_X$ and $\Delta_B'$ of $\Delta_B$ with
$\Delta_B'$ nonsingular, such that the induced map $f':X'\rar B'$ is an
equidimensional toroidal morphism.

If $f^{-1}(U_B)=U_X$ then $(f')^{-1}(U_{B'})=U_{X'}$.
\end{prp}

{\bf Proof.} There exists a subdivision of $\Delta_B$ ``induced'' by
$\Delta_X$. For $x\in\tau\in\Delta_B$ let $S_x$ be the set of cones
$\sigma\in\Delta_X$ such that $\sigma\cap
f_\Delta^{-1}(x)\neq\{0\}$. Since $f_\Delta$ is surjective, $S_x\neq
\emptyset$. The relation $x\sim y \Leftrightarrow S_x = S_y$ for
$x,y\in\tau$ is clearly an equivalence, hence it defines a partition
of $\Delta_B$. The equivalence class of $x$ is
\[ \bigcap_{\sigma\in S_x} f_\Delta(\sigma) \]
which is a convex rational polyhedral subcone of $\tau$. Thus the
partition defines a subdivision $\Delta_B^0$ of $\Delta_B$ such that
$f_\Delta(\sigma)$ for any cone $\sigma\in\Delta_X$ is a union of
cones in $\Delta_B^0$. By the previous lemma there exists a refinement
$\Delta_B^1$ of $\Delta_B^0$ which is a projective subdivision of
$\Delta_B$. Finally, we let $\Delta_B'$ be a nonsingular projective
subdivision of $\Delta_B^1$.

For $f_\Delta$ to map cones of $\Delta_{X}$ to cones of $\Delta_B'$,
the complex $\Delta_{X}$ has to be subdivided. Since the subdivision
$\Delta_B'$ of $\Delta_{B}$ is projective, there exists a good
function $\psi$ on $|\Delta_B'|=|\Delta_B|$. The piecewise linear
function $\psi\circ f_\Delta$ defines a projective subdivision
$\Delta_X'$ of $\Delta_{X}$ whose cones map into cones of
$\Delta_B'$. If $\sigma \in \Delta_X'$ then $f_\Delta(\sigma)$ is a
union of faces of some cone $\tau\in\Delta_B$. Since $f_\Delta$ is
linear on $\sigma$, $f_\Delta(\sigma)$ is convex in $\tau$, hence
$f_\Delta(\sigma)$ is a face of $\tau$.

If we assume from the beginning that $f^{-1}(U_B)=U_X$, that means
$f_\Delta^{-1}(0)\cap |\Delta_X| = 0$, then clearly the same is true
for any subdivision $|\Delta_X'| = |\Delta_X|$, hence
$(f')^{-1}(U_{B'})=U_{X'}$.  \qed

\section{Kawamata's trick and reduced fibers}\label{reduced-fibers}

\subsection{Statement of result} The goal in this section is to find a
finite base change, after which
all the fibers in the resulting morphism are reduced. The base change
we perform will not necessarily be toroidal, but the morphism
after base change will still be toroidal.

\begin{prp}\label{prop-red}
Let $U_X\subset X$ and $U_B \subset B$ be projective toroidal
embeddings, and assume that $B$ is nonsingular. Let $f: X\rar B$ be a
surjective equidimensional toroidal morphism with
$f^{-1}(U_B)=U_X$. Then there exists a finite surjective morphism $p:
B'\rar B$ so that, denoting by $X'$ the normalization of $X\times_B
B'$, we have that $B'$ and $X'$ are toroidal embeddings, the projection $f':
X'\rar
B'$ is an equidimensional toroidal morphism with reduced fibers, and
$(f')^{-1}(U_{B'})=U_{X'}$.
\end{prp}

The construction of $X'$ and $B'$ is more explicitly given in
Proposition~\ref{prop-compl}, where  the polyhedral complexes
of $X'$ and $B'$ are also described.

\subsection{The toric pictures} To start, we characterize equidimensional
toroidal morphisms with reduced fibers in terms of polyhedral
complexes.

\begin{lem} Let $f:X\to B$ be an equidimensional toroidal morphism,
$f_\Delta:\Delta_X\to \Delta_B$ the associated morphism of polyhedral
complexes. Then $f$ has reduced fibers if and only if for any cone
$\sigma \in \Delta_X$, with image $\tau\in \Delta_B$, we have
$f_\Delta(N_\sigma\cap\sigma) = N_\tau\cap\tau$. That is, the image of
the lattice in any cone of $X$ is the lattice in the image cone.
\end{lem}

{\bf Proof.} It suffices to consider the toric morphism of local
models $f: X_\sigma \rar X_\tau$ and the fiber over a point $t\in
X_\tau$ lying in the closed orbit of $X_\tau$. If the orbit of $t$ is
not $\{t\}$ then $X_\tau$ is a product $X_\tau = X_{\tau'}\times
\bfg_m^q$ for some $q>0$. Without loss of generality we may then
replace $X_\tau$ by $X_{\tau'}$, and replace $f$ by $p\circ f$ where $p:
X_\tau\rar X_{\tau'}$ is the projection. Indeed, the fiber $(p\circ
f)^{-1}(p(t))$ is isomorphic to $f^{-1}(t)\times \bfg_m^q$, and
$p_\Delta$ gives an isomorphism $p_\Delta: N_\tau\cap\tau \cong
N_{\tau'}\cap\tau'$. Thus we may assume that $\{t\}$ is the unique
closed orbit of $X_\tau$. The ideal of $f^{-1}(t)$ is generated by
$k[f_\Delta^{\vee}(\tau^\vee\cap(M_\tau\setmin \{0\}))] \subset
k[\sigma^\vee\cap M_\sigma]$. The fiber is reduced if and only if the
image $f_\Delta^{\vee}(\tau^\vee\cap(M_\tau\setmin \{0\}))$ is
saturated in $\sigma^\vee\cap M_\sigma$. This happens if and only if
$f_\Delta(N_\sigma\cap\sigma) = N_\tau\cap\tau$. \qed

When $X_\tau$ is nonsingular, the condition of the lemma is that
primitive points of $\sigma$ map to primitive points of $\tau$.

\begin{lem}\label{lem-tor-prod} Let $X_\sigma\rar X_\tau$ be a toric
morphism with $X_\tau$ nonsingular. Let $X_{\tau'}$ be a toric variety
given by $\tau'=\tau$ and $N_{\tau'}\subset N_\tau$ a sublattice of
finite index. Then every irreducible component of the normalization
of $X_\sigma \times_{X_\tau} X_{\tau'}$ is a toric variety
$X_{\sigma'}$ given by the cone $\sigma'=\sigma$ and integral lattice
$N_{\sigma'} = N_\sigma \cap f_\Delta^{-1}(N_{\tau'})$.
\end{lem}

{\bf Proof.}  The ring of regular functions of $X_{\tau'}
\times_{X_\tau} X_\sigma$ is
\[ \co_{X_{\tau'}}\otimes_{\co_{X_\tau}} \co_{X_\sigma} =
k[(\tau')^{\vee}\cap M_{\tau'}]
\otimes_{k[\tau^{\vee}\cap M_\tau]} k[\sigma^{\vee}\cap M_\sigma] =
k[\pi] \]
where $\pi$ is the pushout of $j: \tau^{\vee}\cap M_\tau \to (\tau')^{\vee}\cap
M_{\tau'} $ and $f^\vee_\Delta: \tau^{\vee}\cap M_\tau \to \sigma^{\vee}\cap
M_\sigma$:

\[ \pi = ((\tau')^{\vee}\cap M_{\tau'}) \times (\sigma^{\vee}\cap
M_\sigma) /\sim.\]

Here $\sim$ is  the equivalence relation generated by:

$(v_1,w_1)\sim
(v_2,w_2)$
whenever  there exists $u\in \tau^{\vee}\cap M_\tau$ such that
$ (v_1,w_1) = (v_2,w_2) \pm (u,-f^\vee_\Delta(u)). $

Let $M$ be the abelian group
\[ M = M_{\tau'}\times M_\sigma /((u,-f^\vee_\Delta(u))) |
u\in \tau^{\vee}\cap M_\tau). \]
We will show that the semigroup homomorphism $\iota: \pi \rar M$ is injective.
 Suppose that $\iota(v_1,w_1) = \iota(v_2,w_2)$, where $v_1,v_2 \in
(\tau')^{\vee}\cap M_{\tau'} $ and $w_1,w_2 \in \sigma^{\vee}\cap
M_\sigma$. Say $  (v_1,w_1) = (v_2,w_2) + (u,-f^\vee_\Delta(u))$ for some $u\in
M_\tau$. Let $t_1,\ldots,t_m \in
\tau^\vee$ be a basis of $M_\tau$, so that in this basis $\tau^\vee$ is given
by the inequalities $t_i\geq 0$, for $i=1,\ldots,k$. Write $u=\sum_i
\alpha_it_1$ with $\alpha_i\in \bfz$. Collecting positive and negative terms,
we get $u = u_1-u_2$ with $u_1,u_2\in \tau^{\vee}\cap M_\tau.$ Writing
$v_1,v_2$ in terms of $t_i$ and expanding the equation $v_1=v_2+u_1-u_2$, we
see that $v_2-u_2 \in \tau^{\vee}\cap M_{\tau'}$. Let $(v_3,w_3) =
(v_2-u_2,w_2+f^\vee_\Delta(u_2))$. Clearly
$$(v_1,w_1)\sim (v_3,w_3)\sim (v_2,w_2)$$
 and
$\iota$ is injective.

The integral closure of $k[\pi]$ in $k[M]$ is the semigroup algebra
$k[\tilde{\pi}]$ where $\tilde{\pi}$ is the saturation of $\pi$ in
$M$. Write $M=F\oplus T$ where $F$ is the free part and $T$ is
torsion. Then $T\subset\tilde{\pi}$ because $m T=0$ for some
$m>0$. For any $f+t \in \tilde{\pi}$ where $f\in F$ and $t\in T$, we
have $-t\in T\subset \tilde{\pi}$, hence $f\in \tilde{\pi}$. Thus
\[ \tilde{\pi} \cong (\tilde{\pi}\cap F) \oplus T.\]
If $t_i$ are generators of $T$ of order $m_i$, and if $x_i$ is the
image of $t_i$ in $k[\tilde{\pi}]$, then
\[ k[\tilde{\pi}] \cong k[\tilde{\pi}\cap F][\ldots,x_i,\ldots] / (x_i
^{m_i}=1) \] Thus the normalization of $\spec {k[\tilde{\pi}]}$ has
$|T|$ components with $x_i = \zeta_i$ where $\zeta_i$ are $m_i$'th
roots of unity, each component isomorphic to $\spec {k[\tilde{\pi}\cap
F]}$.

Next we show that $F$ can be embedded in $M_\sigma\otimes\bfr$ so that
$M_\sigma\subset F$. The image of the homomorphism $\phi: M\rar
M_\sigma\otimes\bfr$ defined by $\phi(v,w)=f_\Delta^\vee(v)+w$ for
$v\in M_{\tau'}, w\in M_\sigma$ contains $M_\sigma$. Since $M$ and
$M_\sigma$ have the same rank, $\phi$ embeds $F$ in
$M_\sigma\otimes\bfr$ as a lattice of full rank containing
$M_\sigma$. As $f_\Delta^\vee$ takes $\tau^\vee$ into $\sigma^\vee$,
we see that $\phi$ maps $\tilde{\pi}\cap F$ into $\sigma^\vee$ so that
the image contains $\sigma^{\vee}\cap M_\sigma$. Therefore, $\spec
{k[\tilde{\pi}\cap F]}$ is a toric variety defined by the cone
$\sigma'=\sigma$ and integral lattice $N_{\sigma'} = F^\vee$. To
determine $F^{\vee}$, first, we have $M_\sigma\subset F$, hence
$F^{\vee}\subset N_\sigma$; second, $f_\Delta^{\vee}(M_{\tau'})\subset
F$ implies that $f_\Delta(F^{\vee})\subset N_{\tau'}$. Conversely,
$N_\sigma\cap f_\Delta^{-1}(N_{\tau'}) \subset F^{\vee}$, so the two
are equal.  \qed

To get a toric morphism with reduced fibers, one can take a base
change $X_{\tau'} \rar X_\tau$, where $\tau'=\tau$ and
$N_{\tau'}\subset N_{\tau}$ is a sublattice of finite index. By
Lemma~\ref{lem-tor-prod} every component $X_{\sigma'}$ of the
normalization of $X_\sigma \times_{X_\tau} X_{\tau'}$ is a toric
variety defined by the cone $\sigma' = \sigma$ and integral lattice
$N_{\sigma'} = N_\sigma \cap f_\Delta^{-1}(N_{\tau'})$. By a judicious
choice of $N_{\tau'}$ the fibers of $X_{\sigma'}\rar X_{\tau'}$ are
reduced.

\subsection{Kawamata's covering}
To perform a similar base change in the toroidal case we need a
toroidal morphism $B'\rar B$ which ramifies over a divisor $D$ with a
certain index.

\begin{dfn} \label{def-cyclic-cover} Let $L$ be a Cartier divisor on $B$, and
let $D\in|m L|$. Choose a
rational section $s_L$ of $ \co_B(L)$ defining $L$ and a regular section
$s_D$ of  $\co_B(mL)$ defining $D$. Consider the rational function $\phi =
s_D/s_L^m$ on $B$. Then the field $K(B)(\sqrt[m]{\phi})$
depends only on $D$ and $\co(L)$. The normalization of $B$ in
$K(B)(\sqrt[m]{\phi})$ is called {\bf the cyclic cover ramified along $D$
with index $m$. }
\end{dfn}

\begin{rem} Another way to define the cyclic cover, is as the normalization of
$${\cal S}pec_B\sym^\bullet\bigl(\co_B(L)^\vee \bigr)/(\{f- s(f)\}_f)
$$
where we view $s$ as a morphism $s:\co_B(kL)\to\co_B((k+m)L)$.
\end{rem}
When $B$ is nonsingular and $D$ a divisor of normal crossings so that
$s_D=x_1\cdots x_l$ for some local parameters $x_1,\ldots,x_l$, then
the cyclic cover has a local equation
\[ z^m = x_1\cdots x_l. \]
It is nonsingular if and only if $l=1$.

Let $U_B\subset B$ be a nonsingular projective (strict) toroidal
embedding. Then  $B\setmin U_B = \sum D_i$ is a strict divisor of
normal crossings. Consider the data $(D_i, m_i)$ where $D_i$ are the
irreducible components of $B\setmin U_B$ and $m_i$ are positive
integers, $i=1,\ldots,m$.

 A {\it Kawamata covering
package} consists of $(D_i, m_i, H_{ij})$ with the following properties:
\begin{enumerate}
\item $H_{ij}$ are effective reduced nonsingular divisors on $B$, for
$i=1,\ldots,m$, $j=1,\ldots,\dim B$.
\item $\sum_i D_i + \sum_{i,j} H_{ij}$ is a reduced divisor of normal
crossings (in particular, $H_{ij}$ are distinct).
\item $D_i+ H_{ij} \in |m_i L_i|$ for some Cartier divisor $L_i$ for all
$i,j$.
\end{enumerate}
To find $H_{ij}$, let $M$ be an ample divisor. Take a multiple of $M$
if necessary so that $m_i M-D_i$ is very ample for all $i$, and choose
$H_{ij}$ general members in $|m_i M-D_i|$.

Now let $B_{ij}$ be the cyclic cover ramified along $D_i+H_{ij}$ with
index $m_i$; let $B'$ be the normalization of
\[ B_{1,1} \times_B \ldots \times_B B_{m,\dimension{B}}\]
and let $p:B'\rar B$ be the projection.

\begin{lem} [Kawamata] The variety $B'$ is nonsingular, ramified with
index $m_i$ along $D_i+H_{ij}$. The reduction of the inverse image
$p^*(\sum_i D_i + \sum_{i,j} H_{ij})_{\red}$ is a divisor of normal crossings.
\end{lem}

{\bf Proof.} Let $x_i$ be local equations of $D_i$, and let $y_{ij}$
be local equations of $H_{ij}$ at $b\in B$. Then $B_{ij}$ is locally
given by the equation
\[ z_{ij}^{m_i} = x_i y_{ij} \]
It suffices to prove that the normalization of $\times_j B_{ij}$ is
nonsingular for all $i$.  If $b\in D_i$
then since $\bigcup_j H_{ij} \cap D_i = \emptyset$, say $b\notin
H_{i,0}$, and $y_{i,0}$ is a unit in $\co_{B,b}$. The normalization of
the product $\times_j B_{ij}$ is locally defined by the equations
\begin{eqnarray*}
z_{i,1}^{m_i} &=& x_i \\ (\frac{z_{i,j}}{z_{i,1}})^{m_i} &=& y_{i,j}
\hspace{.5in} j=2,\ldots \dim B.
\end{eqnarray*}
Since $x_i, y_{i,2},\ldots,y_{i,\dim B}$ are either units or local
parameters in $\co_{B,b}$, it follows that the normalization of the
product $\times_j B_{ij}$ is nonsingular at $b$. A similar situation happens
when $b\notin D_i$, since  then $x_i$ is a unit in
$\co_{B,b}$. Thus, $B'$ is
nonsingular and ramified with index $m_i$ along $D_i$ and
$H_{ij}$. Replacing $B$ by $D_i$ or $H_{ij}$, we get that $p^*(\sum_i
D_i + \sum_{i,j} H_{ij})_{\red}$ is a divisor of normal
crossings. \qed

Let $\tilde{U}_B = B\setmin (\bigcup_{i,j} H_{ij}\cup D_i)$, and
$\tilde{U}_{B'} = p^{-1} (\tilde{U}_B)$. Both $\tilde{U}_B \subset B$
and $ \tilde{U}_{B'} \subset B'$ are toroidal embeddings because $B,
B'$ are nonsingular and the divisors $B\setmin\tilde{U}_B$,
$B'\setmin\tilde{U}_{B'}$ cross normally. From the construction and local
equations we
also see that $p$ is a toroidal morphism with respect to this structure. If
$b'\in B'$, $b=p(b')\in
B$, then there exist local parameters $\{x_l\}$ at $b$ and $\{z_l\}$
at $b'$ such that $p$ is defined by $x_l = z_l^{a_l}$ where $a_l =
m_i$ if the divisor defined by $x_l$ is $D_i$ or $H_{ij}$, otherwise
$a_l =1$. We denote the polyhedral complexes of $\tilde{U}_B \subset
B$ and $ \tilde{U}_{B'} \subset B'$ by $\tilde{\Delta}_B$ and
$\tilde{\Delta}_{B'}$, respectively. So, every cone $\tau' \in
\tilde{\Delta}_{B'}$ is mapped homeomorphically to a cone $\tau \in
\tilde{\Delta}_{B}$ by $p_\Delta$. If $N_\tau$ has basis
$u_1,\ldots,u_m$ then $N_{\tau'}$ has basis $a_1 u_1, \ldots, a_m u_m$
where $a_l$ are as above.

\subsection{The toroidal picture} We need to see how adding the divisors
$H_{ij}$ to the toroidal structure of $B$ affects the toroidal structure of
$X$.

\begin{lem} Let $f: X\rar B$ be toroidal, $B$ nonsingular, H a generic
hyperplane section of $B$. Let $x\in X$, $b=f(x)\in H$. Then
\begin{itemize}
\item[(i)] there exist local models $(X_\sigma = X_\sigma'\times
\bfg_m, z\times 1)$ at $x$ and $(X_\tau = X_\tau'\times \bfg_m,
y\times 1)$ at $b$ such that $H$ corresponds to the divisor
$X_\tau'\times \{1\}$ and the morphism $f$ is a product of toroidal
morphisms
\[ f=g \times id:  X_\sigma'\times \bfg_m \rar X_\tau'\times \bfg_m \]
\item[(ii)] $U_B\setmin H\subset B$, $U_X\setmin f^{-1}(H)\subset
X$ are toroidal embeddings and $f$ is a toroidal morphism of these
embeddings.
\end{itemize}
\end{lem}

{\bf Proof.} Let $(X_{\sigma},s)$ and $(X_{\tau},t)$ be local models
at $x\in X$ and $b=f(x)\in B$, respectively, and let $f$ also denote
the toric morphism of the local models defined by
\[ f_\Delta: (N_\sigma,\sigma)\rar(N_\tau,\tau) \]
Clearly $U_B\setmin H\subset B$ is toroidal and we may assume that
it has a local model $(X_{<\tau,v>},t')$, where $<\tau,v>$ is the cone
spanned by $\tau$ and some indivisible $v\in N_{\tau}$, and where $t'$
lies in the unique closed orbit of $X_{<\tau,v>}$. Write $N_{\tau,2}$
for the saturated sublattice of $N_{\tau}$ generated by $v$, and
choose a splitting
\begin{eqnarray*}
N_{\tau} &\cong& N_{\tau,1} \oplus N_{\tau,2} \\ X_{\tau} &\cong&
X_{\tau}' \times \bfg_m \\ X_{<\tau,v>} &\cong& X_{\tau}' \times
\bfa^1
\end{eqnarray*}

Since $f$ is dominant, $\tilde{N}_{\tau,2} = f_\Delta(f_\Delta^{-1}(
N_{\tau,2}))$ is a sublattice of finite index of $N_{\tau,2}$. The
inclusion $\tilde{N}_{\tau,2} \subset N_{\tau,2}$ corresponds to an
\'etale cover of $\bfg_m$. Since the completed local rings are
isomorphic, we may replace the local model by a local model in the
\'etale cover and assume that $f_\Delta^{-1}(N_{\tau,2})$ surjects
onto $N_{\tau,2}$.

Let $N_{\sigma,1}=f_\Delta^{-1}(N_{\tau,1}) \subset N_{\sigma}$ and
let $N_{\sigma,2} \subset N_{\sigma}$ be generated by some $u\in
N_{\sigma}$ such that $f_\Delta(u)=v$. For any $w\in N_{\sigma}$,
$f_\Delta(w) = v_1+v_2 \in N_{\tau,1} \oplus N_{\tau,2}$, there exists
$m\in \bfz$ so that $f_\Delta(m u)=v_2$. Hence $w-m u \in
N_{\sigma,1}$, and $N_{\sigma}=N_{\sigma,1}+N_{\sigma,2}$. Since
$N_{\sigma,1}\cap N_{\sigma,2} = \{0\}$ the sum is direct,
\begin{eqnarray*}
 N_{\sigma} &\cong& N_{\sigma,1} \oplus N_{\sigma,2} \\ X_{\sigma}
 &\cong& X_{\sigma}' \times \bfg_m
\end{eqnarray*}
and $f_\Delta$ maps $N_{\sigma,1}$ to $N_{\tau,1}$ and $N_{\sigma,2}$
to $N_{\tau,2}$. Thus the toric morphism $f$ is a product $f=g\times
h$ where $g:X_{\sigma}'\rar X_{\tau}'$ and $h=id:\bfg_m \rar \bfg_m$
are the toric morphisms induced by the restriction of $f_\Delta$ to
$N_{\sigma,1}$ and $N_{\sigma,2}$, respectively. Since $f_\Delta$ maps
$<\sigma,u>$ to $<\tau,v>$, $f$ extends to a toric morphism
\[ f: X_{<\sigma,v>} \cong X_{\sigma}' \times \bfa^1 @>{g\times id}>>
X_{\tau}' \times \bfa^1 \cong X_{<\tau,v>}\]

For $t=(y,1)\in X_{\tau}' \times \bfg_m \subset X_{\tau}' \times
\bfa^1$ the complete local rings $\hat{\co}_{X_{\tau}' \times
\bfg_m,(y,1)}$ and $\hat{\co}_{X_{\tau}' \times \bfa^1,(y,0)}$ are
isomorphic via translation. The closures of codimension 1 orbits
through $(y,0)$, are those through $(y,1)$ plus $X_{\tau}' \times
\{0\}$ which formally corresponds to $H$. Similarly, for $s=(z,1)\in
X_{\sigma}' \times \bfg_m \subset X_{\sigma}' \times \bfa^1$ the ring
$\hat{\co}_{X_{\sigma}' \times \bfg_m,(z,1)}$ is isomorphic to
$\hat{\co}_{X_{\sigma}' \times \bfa^1,(z,0)}$ via translation, and the
closures of codimension 1 orbits through $(z,0)$, are those through
$(z,1)$ plus $X_{\sigma}' \times \{0\} = f^{*}(X_{\tau}' \times
\{0\})$ which formally corresponds to $f^{*}(H)$.

Thus $(X_{<\sigma,u>},(z,0))$ and $(X_{<\tau,v>},(y,0))$ are local
models at $x\in X$ and $b\in B$, respectively, and the morphism
\[ f: (U_X\setmin f^{-1}(H)\subset X)\rar(U_B\setmin H\subset B)\]
is toroidal.  \qed

It follows from the lemma that $f:(\tilde{U}_B\subset B) \rar
(\tilde{U}_X\subset X)$ is toroidal, where $\tilde{U}_X =
f^{-1}(\tilde{U}_B)$. By Lemma~\ref{lem-tor-prod}, the normalization
$X'$ of $X\times_B B'$ is toroidal. Let $f': X'\rar B'$, $p':X'\rar X$
be the (toroidal) projections, and let $\tilde{\Delta}_X$,
$\tilde{\Delta}_{X'}$ be the polyhedral complexes of $X, X'$.
\[
\begin{CD}
X' @>{p'}>> X\\ @V{f'}VV @VV{f}V\\ B' @>{p}>> B
\end{CD}
\]
Then $p_\Delta'$ maps a cone $\sigma' \in \tilde{\Delta}_{X'}$
homeomorphically to a cone $\sigma \in \tilde{\Delta}_{X}$; the
integral lattice of $\sigma'$ can then be identified with a sublattice
$N_{\sigma'} \subset N_\sigma$.

The following lemma shows that the added divisors $p^{-1}(H_{ij})$ and
$(p\circ f')^{-1}(H_{ij})$ can be removed from the toroidal structures
of $B'$ and $X'$ so that $f'$ (but not $p$ or $p'$) remains toroidal.

\begin{lem}
$p^{-1}(U_B) \subset B'$ and $(p\circ f')^{-1}(U_B) \subset X'$ are
toroidal embeddings and $f'$ is a toroidal morphism of these
embeddings.
\end{lem}

We show how to remove one irreducible divisor $H=H_{ij}$ for some
$i,j$. Since the question is local, choose local models $X_{\sigma}$,
$X_{\tau}$ and $X_{\tau'}$ of $X$, $B$, and $B'$ so that both $f$ and
$p$ are products
\[ f: X_{\sigma} \cong X_{\sigma}' \times \bfa^1 @>{g\times id}>>
X_{\tau}' \times \bfa^1 \cong X_{\tau}\]
\[ p: X_{\tau'} \cong X_{\rho}' \times \bfa^1 @>{q\times r}>>
X_{\tau}'\times \bfa^1  \cong X_{\tau} \]
where all morphisms are toric and $H$ corresponds to $X_{\tau}' \times
\{0\}$. Since $p$ ramifies along $H$, the morphism
$r:\bfa^1\rar\bfa^1$ has degree $m\geq 2$.

The fiber product is then
\begin{eqnarray*} X'' & = & (X_{\sigma}' \times \bfa^1)\times_{X_{\tau}'
\times \bfa^1} (X_{\tau'}' \times \bfa^1) \\
 &  \cong & (X_{\sigma}' \times_{X_{\tau}'}
X_{\tau'}')\times (\bfa^1 \times_{\bfa^1} \bfa^1) \\
& \cong & (X_{\sigma}'
\times_{X_{\tau}'} X_{\tau'}')\times \bfa^1
\end{eqnarray*}
and the projection
$X''\rar X_{\tau'}$ is the product of the projection $X_{\sigma}'
\times_{X_{\tau}'} X_{\tau'}' \rar X_{\tau'}'$ and $id:
\bfa^1\rar\bfa^1$. The divisor $(p\circ f')^{*}(H)$ in $X''$ is
$X_{\sigma}' \times_{X_{\tau}'} X_{\tau'}' \times \{0\}$ and by the
same translation argument as above we can remove $p^{*}(H)$ and
$(p\circ f')^{*}(H)$ from the toroidal structures of $B'$ and $X''$,
respectively. \qed

Let $\Delta_{X'}, \Delta_{B'}$ be the polyhedral complexes of the
embeddings $p^{-1}(U_B) \subset B'$ and $(p\circ f')^{-1}(U_B) \subset
X'$. Removing the divisors $p^{-1}(H_{ij})$ and $(p\circ
f')^{-1}(H_{ij})$ means removing the corresponding edges (and
everything attached to them) from the polyhedral complexes
$\tilde{\Delta}_{B'}$ and $\tilde{\Delta}_{X'}$. As $p_\Delta$ and
$p_\Delta'$ map cones of $\tilde{\Delta}_{B'}$ and
$\tilde{\Delta}_{X'}$ homeomorphically to cones of
$\tilde{\Delta}_{B}$ and $\tilde{\Delta}_{X}$, restrictions of
$p_\Delta$ and $p_\Delta'$ map cones of ${\Delta}_{B'}$ and
${\Delta}_{X'}$ homeomorphically to cones of ${\Delta}_{B}$ and
${\Delta}_{X}$. We summarize the previous constructions in the
following proposition.

\begin{prp}\label{prop-compl}
With the assumptions of Proposition~\ref{prop-red}, let $u_i$ be the
primitive points of $\Delta_B$, and let $m_i>0$ for
$i=1,\ldots,l$. There exists a finite covering $p:B'\rar B$ so that, if
$X'$ is the normalization of $X\times_B B'$ and $f':X'\rar B'$,
$p':X'\rar X$ the two projections, we have
\begin{itemize}
\item[(i)] $U_{X'}\subset X'$ and $U_{B'}\subset B'$ are toroidal
embeddings with polyhedral complexes $\Delta_{X'}$ and $\Delta_{B'}$;
moreover, $f'$ is a toroidal morphism of these embeddings.
\item[(ii)] There exist morphisms
$$p_\Delta:\Delta_{B'}\rar \Delta_{B}$$
and $$p_\Delta':
\Delta_{X'}\rar \Delta_{X},$$
such that
\begin{list}{$\circ$}{}
\item for any cone $\tau'\in \Delta_{B'}$, the morphism $p_\Delta$ maps
$\tau'$ isomorphically to  $\tau$, and  identifies $N_{\tau'}$ with the
sublattice of $N_{\tau}$ generated by $m_i u_i$;
\item  for any cone $\sigma'\in \Delta_{X'}$, $p_\Delta'$ maps
$\sigma'$ isomorphically to $\sigma$, and identifies
$N_{\sigma'}$ with the sublattice $N_\sigma \cap f_\Delta^{-1}
(N_{\tau'})$ of $N_{\sigma}$
\item The following diagrams commute:

\begin{minipage}{2.5in}
\[
\begin{CD}
\sigma' @>{\cong}>p_\Delta'> \sigma\\ @V{f_\Delta'}VV @VV{f_\Delta}V\\
\tau' @>{\cong}>p_\Delta> \tau
\end{CD}
\]
\end{minipage}
\begin{minipage}{2.5in}
\[
\begin{CD}
N_{\sigma'} @>{\subset}>p_\Delta'> N_\sigma\\ @V{f_\Delta'}VV
@VV{f_\Delta}V\\ N_{\tau'} @>{\subset}>p_\Delta> N_{\tau}
\end{CD}
\]
\end{minipage} \\
\end{list}
\end{itemize}
\end{prp} \qed

{\bf Proof of \ref{prop-red}} Let $\Delta_X$, $\Delta_B$ be the
polyhedral complexes of $X$, $B$, and $f_\Delta:
\Delta_X\rar\Delta_B$.  By equidimensionality, $f_\Delta$ maps
$\sigma^{(1)}$ to $\tau^{(1)}$. If $u_i$ are the primitive
points  of $\Delta_B$, let $v_{ij}$ be the
primitive points of $\Delta_X$ such that $f_\Delta(v_{ij}) = m_{ij}
u_i$ for some $m_{ij} >0$. Set $m_i = \lcm_j \{m_{ij}\}$. Then for all
$i,j$, some multiple of $v_{ij}$ maps to $m_i u_i$. Now we use the
covering data $(D_i, m_i)$ to define the toroidal morphism $f':X'\rar
B'$ as in the previous proposition. Let $\sigma' \in \Delta_{X'}$,
$f_\Delta'(\sigma') = \tau' \in \Delta_{B'}$. There exist
$\sigma\in\Delta_X$, $\tau\in\Delta_B$ such that $\sigma'\cong\sigma$,
$\tau'\cong\tau$. The integral lattice $N_{\tau'}$ is generated by
$m_i u_i$ with $u_i$ lying on the edges of $\tau$, and
$N_{\sigma'}=N_\sigma \cap f_\Delta^{-1} (N_{\tau'})$. For any edge of
$\sigma'$ spanned by $v_{ij}$, some multiple $a_{ij} v_{ij}$ maps to
$m_i u_i$. Hence $a_{ij} v_{ij} \in N_{\sigma'}$ is primitive and maps
to the primitive point $m_i u_i$. This proves that $f'$ has reduced
fibers. \qed

\section{Mildness of the morphism}\label{mildness}

Recall that a normal variety $X$ is said to have Gorenstein
singularities if it is Cohen-Macaulay and has an invertible dualizing
sheaf $\omega_X$. It is said to have rational singularities, if for  a
resolution of singularities  $r:X'\to X$ we have $r_*\omega_{X'} =
\omega_X$. Every toric variety is Cohen-Macaulay with rational
singularities (\cite{te} Theorem~14, p. 52). The dualizing sheaf of an affine
toric variety $X_\tau$ is the coherent sheaf $\omega_{X_\tau} =
\co(-\sum D_i)$ where $D_i$ are closures of the codimension 1
orbits. This sheaf is invertible if and only if there exists an element of $M$
(namely a linear function on $\tau$ taking integer values on $N$) such that
$\psi(v) = -1$ for all  $v\in \tau^{(1)}$.

\begin{lem}
Let $X_\sigma$ and $X_\tau$ be affine toric varieties with $X_\tau$
nonsingular. Let $f: X_\sigma \rar X_\tau$ be an equidimensional toric
morphism without horizontal divisors, having only reduced fibers. Then
$X_\sigma$ has rational Gorenstein singularities.
\end{lem}

{\bf Proof.}  Since $X_\tau$ is nonsingular, it has rational
Gorenstein singularities. Let $\psi:\tau\rar\bfr$ be the linear
interpolation of $\psi(u_i)=-1$ for every primitive point $u_i$ of
$\tau$. Then $\psi\circ f_\Delta:\sigma\rar\bfr$ is the required
function. It is linear because $\psi$ and $f_\Delta$ are; since $f$ is
equidimensional and has reduced fibers, primitive points in $\sigma$
map to primitive points in $\tau$, hence $\psi\circ f_\Delta$ takes
the value $-1$ on every primitive point of $\sigma$.  \qed

\begin{lem}\label{lem-rational-Gorenstein}
Let $U_X\subset X$ and $U_B \subset B$ be toroidal embeddings, and
assume that B is nonsingular. Let $f: X\rar B$ be an equidimensional
toroidal morphism, without horizontal divisors, and with reduced fibers. Then
$X$ has rational Gorenstein singularities.
\end{lem}

{\bf Proof.}  Having rational Gorenstein singularities is a local
analytic property. Since all local models have rational Gorenstein
singularities, so does $X$.  \qed

To show that the morphism is mild, we need to look at the situation
after dominant base changes. We first look at the cases where the base
change is relatively nice:

\begin{lem}\label{lem-log-smoothness}
 Let $f:X\to B$ be as above. Let $g:B'\to B$ be a dominant morphism,
where $B'$ a  nonsingular variety,  and assume that $ g^{-1} (B\setmin U_B)$ is
a normal
crossings divisor. Let $X'$ be the pullback of $X$ to $B'$. Then
$X'\to B'$ admits a toroidal structure relative to $ g^{-1}U_B \subset B'$.
\end{lem}

{\bf Proof.} We use the formalism of logarithmic structures. By \cite{kato1},
\S 8.1,
the morphism $X\to B$ is logarithmically smooth. Also $ g^{-1}U_B \subset B'$
endows $B'$ with a logarithmically regular structure. Moreover $g:B'\to B$ is a
morphism of logarithmic schemes. The variety $X'$ thus inherits a logarithmic
structure. The morphism $X'\to B'$ clearly satisfies the formal lifting
property for logarithmic smoothness (\cite{kato1},
\S 8.1,(i)). It is left to show that $X'$ satisfies condition (S)
(\cite{kato1} \S 1.5). Indeed, since $f$ is equidimensional, the fibers of $f$
are
reduced and $B'$ is normal, it follows that $X'$ is regular in codimension
1. Since $f$ is a Gorenstein morphism and $B$ is Gorenstein, we have that $X'$
is
Gorenstein, in particular it is Cohen Macaulay. It then follows that $X'$ is
normal. Combining this with the assumption that $f^*M_B$ is saturated in $M_X$,
we have that the monoids giving the logarithmic charts
on $X'$ are integral and saturated. Altogether, we have that $X'$ satisfies
condition (S). Therefore $f':X'\to B'$ is logarithmically smooth, and thus
toroidal.

\begin{lem}\label{lem-descent} Let $f:X\to B$ be a flat Gorenstein morphism and
assume that $B$
has
rational Gorenstein singularities. Assume there is a modification $r:B'\to B$
with rational Gorenstein singularities, such that $X'=X\times_B B'$ has
rational singularities. Then $X$ has rational singularities as well.
\end{lem}

{\bf Proof.} Consider the diagram
$$\begin{array}{rcl} X' &\stackrel{h}{\to}  & X \\
               f' \downarrow & &   \downarrow f\\
             B' &\stackrel{r}{\to}  & B.
\end{array}$$

By base change we have $h^*\omega_f = \omega_{f'}$. By assumption we have
$\omega_B = r_*\omega_{B'}$. The flatness of $f$ implies $f^*r_*\omega_{B'} =
h_*{f'}^*\omega_{B'}$. Therefore we have
\begin{eqnarray*}
h_*\omega_{X'} & = & h_*(\omega_{f'} \otimes {f'}^*\omega_{B'} )\\
 & = &  h_*(h^*\omega_f \otimes {f'}^*\omega_{B'}) \\
 & = & \omega_f \otimes h_*{f'}^*\omega_{B'} \\
 & = &  \omega_f \otimes f^*r_*\omega_{B'} \\
 & = &  \omega_f \otimes f^*\omega_B \\
 & = & \omega_X.
\end{eqnarray*}
Thus $X$ has rational singularities. \qed

\begin{prp}
Let $U_X\subset X$ and $U_B \subset B$ be toroidal embeddings and
assume that B is nonsingular. Let $f: X\rar B$ be an equidimensional
 morphism with reduced fibers, which is toroidal  such that  $U_X=f^{-1}U_B$.
 Then $f$ is mild.
\end{prp}

{\bf Proof.} Let $B_1\to B$ be a dominant morphism such that $B_1$ has rational
Gorenstein singularities. We need to show that $X_1 = X\times_BB_1$ has
rational Gorenstein singularities. By lemma \ref{lem-rational-Gorenstein} it
has Gorenstein singularities.

Pick a resolution of singularities $B'\to B_1$ such that the inverse image of
$B\setmin U_B$ in $B'$ is a divisor with normal crossings. By lemma
\ref{lem-log-smoothness} we have that $X' = X\times_BB'$ is toroidal, therefore
$X'$ has rational  singularities. By lemma \ref{lem-descent} we have that $X_1$
has rational singularities as well, which is what we needed. \qed

\section{The Cartier covering}\label{sec-cartier} In section
\ref{combinatorial} below we will
translate the semistable reduction conjecture in purely combinatorial terms. In
order to maximize the flexibility of the combinatorial operations, we need to
generalize Kawamata's trick slightly to accommodate cases where $B$ has
quotient singularities. In such a situation we need a finite cover of $B$ which
is nonsingular, and such that the resulting polyhedral complex is easily
described.  The following statement will suffice for this purpose.

\begin{prp}\label{prop-cartier} Let $U_X\subset X$ and $U_B \subset B$ be
projective toroidal
embeddings, and assume that $B$ has only quotient singularities.
Let $f: X\rar B$ be a surjective mild toroidal morphism,
satisfying $f^{-1}(U_B)=U_X$. Then there exists a finite surjective morphism
$p: B'\rar B$ so that, denoting  $X'=X\times_B B'$, $U_{B'}=p^{-1}U_B$ and
$U_{X'}= X\times_B U_{B'}$ we have
\begin{enumerate}
\item $B'$ and $X'$ are toroidal embeddings with polyhedral complexes
$\Delta_{X'}$ and $\Delta_{B'}$, with $B'$ nonsingular;
\item the projection $f': X'\rar B'$ is a mild toroidal morphism;
\item  There exist morphisms
$$p_\Delta:\Delta_{B'}\rar \Delta_{B}$$
and $$p_\Delta':
\Delta_{X'}\rar \Delta_{X},$$
such that
\begin{list}{$\circ$}{}
\item for any cone $\tau'\in \Delta_{B'}$, the morphism $p_\Delta$ maps $\tau'$
isomorphically to  $\tau$, and  identifies
$N_{\tau'}$ with the sublattice of
$N_{\tau}$ generated by the primitive vectors of $\tau$;
\item  for any cone $\sigma'\in \Delta_{X'}$, $p_\Delta'$ identifies
$N_{\sigma'}$ with the sublattice
$N_\sigma \cap f_\Delta^{-1} (N_{\tau'})$ of $N_{\sigma}$, and maps
$\sigma'$ isomorphically to $\sigma$.
\item The following diagrams commute:

\begin{minipage}{2.5in}
\[
\begin{CD}
\sigma' @>{\cong}>p_\Delta'> \sigma\\ @V{f_\Delta'}VV @VV{f_\Delta}V\\
\tau' @>{\cong}>p_\Delta> \tau
\end{CD}
\]
\end{minipage}
\begin{minipage}{2.5in}
\[
\begin{CD}
N_{\sigma'} @>{\subset}>p_\Delta'> N_\sigma\\ @V{f_\Delta'}VV
@VV{f_\Delta}V\\ N_{\tau'} @>{\subset}>p_\Delta> N_{\tau}
\end{CD}
\]
\end{minipage} \\
\end{list}
\end{enumerate}
\end{prp}

{\bf Proof.} We use a construction analogous to Kawamata's covering. First, let
$B$ be a variety, $E$ an effective Weil divisor on $B$. We make the following
assumptions:
\begin{enumerate}
\item as a scheme, $E$ is integral and normal;
\item the divisor $ mE = D$ is a Cartier;
\item There is a Cartier divisor  $L$ on $B$ such that $\co(D) = \co(mL)$.
\end{enumerate}
As in  Definition \ref{def-cyclic-cover}, we choose a function $\phi$ with
$\operatorname{div}(\phi) = D-mL$, and define the cyclic cover $p:B'\to B$ by
taking $m$-th root of $\phi$. The point is, that the $\bfq$-Cartier divisor
$p^*E$ is in fact Cartier.

Let us see what happens in the toric case.

\begin{lem} Let $X_\tau$ be an affine toric variety with quotient
singularities and $E$ a toric divisor corresponding to a primitive vector $v$.
Write $\tau=\bfr_+v\times \tau_1$. Choose  $m\in \bfn$ such that $D = mE$ is
Cartier. Then $\co(mD)$ is trivial, and the corresponding cyclic cover is
$X_{\tau'}$, where $|\tau'| = \tau$ and $N_{\tau'} = \bfz v\times
(N_\tau\cap\tau_1)$.
\end{lem}

This lemma is easy and left to the reader. If we iterate the lemma with all the
toric divisors, we obtain a nonsingular covering, whose lattice is generated by
the primitive vectors of $\tau$.

Assume that in addition we have a mild morphism $X_\sigma\to X_\tau$.
Lemma \ref{lem-tor-prod} works word-for-word in this case. The only thing which
needs to be changed in the proof is, that since $f_\Delta(N_\sigma) = N_\tau$,
the semigroup homomorphism $i:\pi\to M$ is still injective. (This is not a
serious business - we could replace the product by its reduction anyway.)

In order to define a Kawamata covering package, we need the following Bertini
type lemma (this is a special case of the stratified Bertini Theorem) .

\begin{lem}
Let $U_B\subset B$ be a toroidal embedding. Let $H_0$ be a very ample divisor
on $B$ and let $H$ be a general element of $|H_0|$. Let $U_B' = U_B\setmin
H$. We have
\begin{enumerate}
\item $U_B'\subset B$ is a toroidal embedding;
\item let $b\in H$ and assume that $b$ lies on a stratum $Y\subset B$ of the
toroidal embedding  $U_B\subset B$. Then $H\cup Y$ is nonsingular at $b$;
\item  there is a local model $X_\rho \times \bfg_m^k$ for $U_B\subset B$ at
$b$, where
$Y$ corresponds to the factor $\bfg_m^k$, and $H$ corresponds to  $X_\rho
\times \bfg_m^{k-1}\times\{1\}$.
\item  there is a local model $X_\rho \times \bfg_m^{k-1}\times\bfa^1$ for
$U_B'\subset B$ at
$b$, where $Y$ corresponds to $\bfg_m^k\times\bfa^1$, and $H$ corresponds to
$X_\rho
\times \bfg_m^{k-1}\times\{0\}$.
\end{enumerate}
\end{lem}

{\bf Proof.}  Since $B$ is a stratified space, it has a local product structure
$X_\rho\times Y$ at $b$. We need to adapt this
structure to  the divisor $H$ as in the lemma. Applying Bertini's theorem to
the closure of the stratum $Y$, we have that $H\cup Y$ is nonsingular. Pick a
regular function $u$ near $b$ defining $H$ on a neighborhood of $b$.  Then the
restriction $u_Y$ of $u$ to $Y$ is a regular parameter. We can find a regular
system of parameters on $Y$ at $b$ including $u_Y$ and lift them to $B$. This
lifting gives the desired adaptation of the product structure. \qed

Using this lemma, the proof of Proposition \ref{prop-compl} works word-for-word
here, and yields Proposition \ref{prop-cartier}. \qed

\section{Towards semistable reduction}\label{combinatorial}
\subsection{Combinatorial semistable reduction}
We are going to restate the semistable reduction conjecture (Conjecture
\ref{conj-semistable}) in purely combinatorial terms. First let us define
 semistability combinatorially.
\begin{dfn}\label{dfn-comb-semi} Let $\Delta_X$ and $\Delta_B$
be rational conical polyhedral complexes, and let $f_\Delta:\Delta_X\to
\Delta_B$ be a surjective polyhedral map.
We say that $f_\Delta$ is {\bf weakly semistable} if the following conditions
hold.
\begin{enumerate}
\item We have $f_\Delta^{-1}(0) = \{0\}$.
\item For any cone $\tau\in \Delta_X$, we have $f_\Delta(\tau)\in \Delta_B$.
\item\label{item-reduced} For any cone $\tau\in \Delta_X$, we have
 $f_\Delta(N_\tau) =
 N_{f_\Delta(\tau)}$.
\item $\Delta_B$ is simplicial with index 1.
\end{enumerate}
If also $\Delta_X$ is simplicial with index 1, we say that  $f_\Delta$ is {\bf
semistable}.
\end{dfn}

We now define the operations we are allowed to perform.
\begin{dfn} Let $\Delta$ be a rational conical polyhedral complex. A {\bf
lattice alteration} $\Delta_1\to \Delta$ consists of an integral structure
induced by a consistent choice of sublattices $N^1_\sigma\subset N_\sigma$ for
each cone $\sigma\subset \Delta$. An {\bf alteration} $\Delta_1\to \Delta$ is a
composition $\Delta_1\to \Delta'\to  \Delta$ of  a lattice alteration
$\Delta_1\to \Delta'$  with a subdivision $\Delta'\to \Delta$. The alteration
is {\bf projective} if the corresponding subdivision $\Delta'\to \Delta$ is
projective.
\end{dfn}

Note that the composition of alterations is an alteration. Be warned that the
factorization of a polyhedral alteration above is not analogous to the Stein
factorization in algebraic geometry - it has the opposite order.

\begin{dfn} \begin{enumerate}
\item  Let $\Delta_X\to \Delta_B$ be a polyhedral map of  rational conical
polyhedral complexes. Let $\Delta_B'\to \Delta_B$ be a subdivision. The induced
subdivision $\Delta_X'\to \Delta_X$ is the minimal subdivision admitting a
polyhedral map to $\Delta_B'$. The cones of $\Delta_X'$ are of the form
$\sigma\cap f_\Delta^{-1}(\tau')$ for $\sigma\in \Delta_X$ and $\tau'\in
\Delta_B'$.
\item Let $\Delta_X\to \Delta_B$ be a polyhedral map of  rational conical
polyhedral complexes. Let $\Delta_B^1\to \Delta_B$ be a lattice alteration. The
induced lattice alteration $\Delta_X^1\to \Delta_X$ is the minimal lattice
alteration mapping to $\Delta_B^1$. The lattices of $\Delta_X^1$ are of the
form $N_\sigma\cap f_\Delta^{-1}(N_{\tau^1})$ for $\sigma\in \Delta_X$ and
$\tau^1\in \Delta_B^1$.
\item  Let $\Delta_X\to \Delta_B$ be a polyhedral map of  rational conical
polyhedral complexes. Let $\Delta_B^1\to \Delta_B'\to \Delta_B$ be an
alteration, factored as a lattice alteration followed by a subdivision. The
induced alteration $\Delta_X^1\to \Delta_X$ is the induced lattice alteration
$\Delta_X^1\to \Delta_X'$ of the induced subdivision  $\Delta_X'\to \Delta_X$.
 \end{enumerate}\end{dfn}

Note that an alteration induced by a projective alteration is projective.

We are now ready to state our conjecture.

\begin{conj}\label{conj-comb-semistable}
 Let $f_\Delta:\Delta_X\to \Delta_B$ be a polyhedral map of  rational conical
polyhedral complexes, and assume for simplicity that $f_\Delta^{-1}(0) =
\{0\}$. Then there exists a projective alteration $\Delta_B^1\to \Delta_B$,
with induced
alteration  $\Delta_X^1\to \Delta_X$, and a projective subdivision $\Delta_Y\to
\Delta_X^1$, such that $\Delta_Y \to \Delta_B^1$ is semistable.
\end{conj}

The conjectures are tied together by the following proposition:

\begin{prp}
Conjecture \ref{conj-comb-semistable} implies Conjecture \ref{conj-semistable}.
\end{prp}

{\bf Proof.} Let $X\to B$ be a surjective  morphism  of complex projective
varieties with geometrically integral generic fiber. By Theorem
\ref{th-toroidal-reduction} we may assume the morphism toroidal. Let
$f_\Delta:\Delta_X \to \Delta_B$ be the associated polyhedral map. By
theorem \ref{th-weak-semistable-reduction} we may assume $f_\Delta$ is weakly
semistable. We now invoque Conjecture \ref{conj-comb-semistable}. Let
$\Delta_B^1\to \Delta_B'\to\Delta_B$ be a projective alteration, $\Delta_X^1\to
\Delta_X$ the induced alteration, and $\Delta_Y\to
\Delta_X^1$ a projective subdivision such that $\Delta_Y \to \Delta_B^1$ is
semistable. The subdivision $\Delta_B'\to\Delta_B$ and the induced subdivision
$\Delta_X'\to\Delta_X$ give rise to  modifications
$B'\to B$ and  $X'\to X$ with a lifting $X'\to B'$ (see lemma
\ref{lem-lifting}). Let $\Delta_B^0\to \Delta_B'$ be the lattice alteration
given by the sublattice generated by the primitive vectors, and $\Delta_X^0\to
\Delta_X'$ the induced lattice alteration.  Proposition
\ref{prop-cartier} gives rise to an alteration $\widetilde{B_0} \to B'$, with
pullback $\widetilde{X_0} \to X'$, and a
polyhedral map $\Delta_{\widetilde{B_0}}\to \Delta_B^0$ which is an isomorphism
on each cone. The same holds for $\Delta_{\widetilde{X_0}}\to \Delta_X^0$. The
lattice alteration $\Delta_B^1\to \Delta_B'$, being of index 1, factors through
$\Delta_B^0$. Proposition \ref{prop-compl} gives an alteration
$\widetilde{B_1}\to\widetilde{B_0}$ with pullback
$\widetilde{X_1}\to\widetilde{X_0}$ and a polyhedral map
$\Delta_{\widetilde{B_1}}\to  \Delta_B^1$ which is an
isomorphism
on each cone, and the same holds for $\Delta_{\widetilde{X_1}}\to
\Delta_X^1$. The subdivision $\Delta_Y\to \Delta_X^1$ induces a subdivision
$\widetilde{\Delta_{Y}}\to\Delta_{\widetilde{X_1}}$, giving rise to a toroidal
modification $\widetilde{Y}\to \widetilde{X_1}$. Since $\Delta_Y\to
\Delta_B^1$ is semistable, we have that
$\widetilde{\Delta_{Y}}\to\Delta_{\widetilde{B_1}}$ is semistable as well, and
therefore $\widetilde{Y}\to\widetilde{B_1}$ is semistable, as required. \qed

\subsection{How far can we push the results?} In \cite{ar}, it is shown that if
$\Delta_0\subset \Delta$ is a subcomplex, and $\Delta_0'\to \Delta_0$ is a {\em
projective} triangulation, then there is a projective triangulation $\Delta'\to
\Delta$
extending $\Delta_0'$, without new edges: ${\Delta'}^{(1)} = {\Delta}^{(1)}
\cup {\Delta_0'}^{(1)}$. Applying the main theorem of \cite{te}, Chapter III,
it is shown in \cite{ar} that given a weakly semistable $X\to B$, one can
locally find a finite
map $B_1\to B$ and a modification $X_2\to X\times_BB_1$ such that $X_2\to B_1$
is semistable in codimension 1. If we apply Kawamata's trick, then clearly this
can be done globally. Moreover, $X_2$ has only quotient singularities.

A weakly
semistable morphism is said to be {\bf almost semistable} if it is semistable
in codimension 1 and moreover $X$ has quotient singularities. Thus, in theorem
\ref{th-weak-semistable-reduction} we may replace ``weakly semistable'' by
``almost semistable''. An analogous definition can be made on the polyhedral
side.

It is important to note that an almost semistable morphism is not necessarily
semistable. It is easy to give a polyhedral example: let $\tau = (\bfr_+)^2$ be
endowed with the standard lattice $\bfz^2$, and let $\sigma = (\bfr_+)^4$
be given the lattice generated by  $\bfz^4$ and the vector
$w=(1/2,1/2,1/2,1/2)$.
We have a polyhedral map $\sigma\to \tau$ given by $(a,b,c,d)\mapsto (a+b,
c+d)$. It is easy to see that this is almost semistable, but not semistable,
since $\sigma$ has index 2. Needless to say, a corresponding toroidal example
can be easily constructed as well.

The example we just gave is easy to amend. Indeed, if we subdivide $\tau$ at
its barycenter $(1,1)$, take the induced subdivision of $\sigma$, then its
star subdivision centered at $w$, and extend this to a triangulation using
\cite{ar}, we obtain a semistable map. This can be extended to families of
surfaces in general - the main observation (see \cite{wang} for the case
$\dim(B)=1$) is  that one can  use Pick's theorem and subdivide, with no need
for additional lattice alteration.  We plan to pursue this elsewhere.

One last remark: the second author has shown, that in order to prove semistable
reduction, it is sufficient to produce $B_1$ and $Y$ such that $Y\to B_1$
satisfies all but condition \ref{item-reduced} of the requirements for
semistability in Definition \ref{dfn-comb-semi}. Again, this will be pursued
elsewhere.


\begin{thebibliography}{HHHHHHH}

\bibitem[$\aleph$]{fibered} D. Abramovich, {\em A high fibered power
of a family of varieties of general type dominates a variety of
general type}, Inventionnes Math., to appear. \\ {\tt
alg-geom/9604024}

\bibitem[$\aleph$-dJ]{aj} D. Abramovich and A. J. de Jong, {\em
Smoothness, semistability and toroidal geometry}, J. Alg. Geom., to
appear.\\{\tt alg-geom/9603018}

\bibitem[$\aleph$-R]{ar}  D. Abramovich and J. M. Rojas, {\em Extending
triangulations and semistable reduction}, preprint.

\bibitem[Al]{alex} V. Alexeev, {\em Boundedness and $K\sp 2$ for log
surfaces.} Internat. J. Math. 5 (1994), no. 6, 779-810.

\bibitem[Al1]{alex2} V. Alexeev, {\em Moduli spaces $M_{g,n}(W)$ for
surfaces,} preprint.\\ {\tt alg-geom/9410003}

\bibitem[A-W]{aw} M. Artin and G. Winters, {\em Degenerate fibers and
reduction of curves,} Topology 10 (1971) p. 373--383.

\bibitem[B-M]{bm} E. Bierstone and P. Milman, {\em Canonical
    desingularization in characteristic zero by blowing up the maximum
    strata of a local invariant}, nvent. Math. 128 (1997), no. 2,
207--302.

\bibitem[D]{danilov} V. I. Danilov, {\em The Geometry of Toric
Varieties}, Russian Math. Survey, 33:2 (1978), 97-154.

\bibitem[dJ]{dj} A. J. de Jong, {\em Smoothness, semistability, and
alterations}, Publ. Math. I.H.E.S. {\bf 83}, 1996, pp. 51-93.

\bibitem[dJ2]{dj2} A. J. de Jong, {\em Families of curves and
alterations},  Annales de l`Institute Fourier, to appear.

\bibitem[Kawa]{kawamata} Y. Kawamata, {\em Characterization of Abelian
Varieties,} Compositio mathematica, vol. 43, Fasc. 2, 253-276 (1981).

\bibitem[Kato]{kato}  K. Kato, {\em logarithmic structures of
Fontaine-Illusie},  Algebraic analysis, geometry and number theory: Proceedings
of the JAMI Inaugural Conference, 191--224, Johns Hopkins Univ. Press,
Baltimore, MD, 1989.

\bibitem[Kato1]{kato1}  K. Kato, {\em Toric singularities.} Amer. J. Math. 116
(1994), no. 5, 1073--1099.

\bibitem[KKMS]{te} G. Kempf, F. Knudsen, D. Mumford and
B. Saint-Donat, {\em Toroidal Embeddings I}, Springer, LNM 339, 1973.

\bibitem[Mor]{m} R. Morelli, {\em The birational geometry of toric
varieties}. J. Alg. Geom. {\bf 5} (1996) 751-782.

\bibitem[Vie]{v} E. Viehveg, {\em Quasi-projective moduli of polarized
manifolds}, Ergebnisse der Mathematik und ihrer Grenzgebiete Band
30. Springer, 1995.

\bibitem[Vil]{villa} O. Villamayor,{\em Constructiveness of Hironaka's
resolution.} Ann. Sci. \'Ecole Norm. Sup. (4) 22 (1989), no. 1, 1--32.

\bibitem[Wang]{wang} J.-H. Wang, MIT dissertation.

\bibitem[W{\l}o]{w} J. W{\l}odarczyk, {\em Decomposition of birational toric
maps in blow-ups and  blow-downs.} Trans. Amer. Math. Soc. 349 (1997), no. 1,
373--411.

\end{thebibliography}
\end{document}